\documentclass[letterpaper,11pt]{article}

\usepackage{ucs}
\usepackage[utf8x]{inputenc}
\usepackage{graphicx}
\usepackage{amsfonts}
\usepackage{dsfont}
\usepackage{amssymb}
\usepackage{amsmath}
\usepackage{amsthm}
\usepackage{enumerate}
\usepackage{stmaryrd}
\usepackage{fullpage}
\usepackage{ifthen}
\usepackage{subfigure}
\usepackage{epic}
\usepackage{authblk}
\usepackage{textcomp}
\usepackage[small]{caption}

\usepackage[hypertexnames=false,colorlinks=true,linkcolor=blue,citecolor=blue]{hyperref}
\usepackage[numbers,comma,square,sort&compress]{natbib}
\usepackage[letterpaper,text={7in,9in},centering]{geometry}

\usepackage{color}
\usepackage{titlesec}
\graphicspath{{eps/}{pdf/}}
\newcommand{\bqq}{\begin{equation}}
\newcommand{\eqq}{\end{equation}}
\newcommand{\bqs}{\begin{equation*}}
\newcommand{\eqs}{\end{equation*}}

\newcommand{\R}{\mathbb{R}} 
\newcommand{\Z}{\mathbb{Z}}

\newcommand{\rme}{\mathrm{e}}
\newcommand{\rmi}{\mathrm{i}}
\newcommand{\rmd}{\mathrm{d}}
\newcommand{\rmo}{{\scriptstyle\mathcal{O}}}
\newcommand{\rmO}{\mathcal{O}}
\newcommand{\eps}{\varepsilon}

\numberwithin{equation}{section}

\newenvironment{Hypothesis}[1]%
  {\begin{trivlist}\item[]{\bf Hypothesis #1 }\em}{\end{trivlist}}

\theoremstyle{plain}
\newtheorem{theorem}{Theorem}[section]
\newtheorem{proposition}[theorem]{Proposition}
\newtheorem{lemma}[theorem]{Lemma}
\newtheorem{corollary}[theorem]{Corollary}

\newtheorem{main}[theorem]{Main Result}

\newenvironment{Proof}[1][.]%
 {\begin{trivlist}\item[]\textbf{Proof#1 }}%
 {\hspace*{\fill}$\rule{0.3\baselineskip}{0.35\baselineskip}$\end{trivlist}}

\title{Pinning and Unpinning in Nonlocal Systems \footnote{This research was conducted during Summer 2014 in the REU: Complex Systems at the University of Minnesota Department of Mathematics, funded by the National Science Foundation (DMS-1311414) and (DMS-1311740); see (\href{http://math.umn.edu/~gfaye/reu.html}{\textcolor{blue}{http://math.umn.edu/$\sim$gfaye/reu.html}}).}}
\author[1]{Taylor Anderson}
\author[2]{Gr\'egory Faye}
\author[2]{Arnd Scheel\footnote {Corresponding author email: \textsc{scheel@umn.edu}}}
\author[3]{David Stauffer}
\affil[1]{\small Mount Holyoke College,
Department of Mathematics and Statistics, 50 College St,
South Hadley, MA 01075, USA}
\affil[2]{\small University of Minnesota,
School of Mathematics, 206 Church St SE, Minneapolis, MN 55414, USA}
\affil[4] {Cornell University, 
Department of Mathematics, 
310 Malott Hall, 
Ithaca, NY 14853, USA}

\begin{document}
\maketitle

\begin{abstract}
\noindent We investigate pinning regions and unpinning asymptotics in nonlocal equations. We show that phenomena are related to but different from pinning in discrete and inhomogeneous media. We establish unpinning asymptotics using geometric singular perturbation theory in several examples. We also present numerical evidence for the dependence of unpinning asymptotics on regularity of the nonlocal convolution kernel. 
\end{abstract}

{\noindent \bf Keywords:} front pinning, singular perturbations, traveling waves, nonlocal coupling\\

\section{Introduction}

Relaxation to the energy minimum in spatially extended systems is often mediated by the propagation of interfaces, separating globally and locally minimizing states. The speed of propagation of such interfaces gives crucial information on time scales for relaxation. In the simplest, typical scenario, the interface motion is driven by the energy difference between local and global minimizers, yielding an effective force on the interface. The speed of the front is then proportional to this effective force. Stationary interfaces correspond to the situation where the states on either side of the interface have equal energy. In formulas, the speed $c$ depends in a smooth and monotone fashion on the energy difference $\mu$, $c=c(\mu)$, $c'(\mu)>0$, so that for $\mu$ small, 
\begin{equation}\label{e:1/1}
c\sim \mu.
\end{equation}
A prototypical example for this scenario is the Nagumo equation,
\begin{equation}\label{e:nag}
u_t=u_{xx}+f(u),\qquad f(u)=u(1-u)(u-a),
\end{equation}
where $\mu= a-\frac{1}{2}$, and $c=\sqrt{2}\mu$ is linear in $\mu$. Similar relations hold for general bistable nonlinearities. 

It has been well known that this simple picture fails in many important contexts. In particular, the speed of interfaces may vanish for sufficiently small yet non-zero energy differentials. Such propagation failure is usually referred to as \emph{pinning}, alluding to a simple scenario where energies depend on space $x$. In the following, we briefly describe this scenario from several view points. Our contribution in this paper can be understood as providing a new, different view point on pinning, with fundamentally different characteristic expansions and analytic tools.

\subsection{Pinning and inhomogeneities}

\paragraph{Periodic media.}

The possibly simplest case where pinning is observed are spatially periodic media, such as 
\begin{equation}\label{e:per}
u_t=u_{xx}+u(1-u)(u-a+\varepsilon \sin(x)).
\end{equation}
The parameter $a$ still detunes the relative energy of the equilibria $u=0$ and $u=1$. Since, however, this relative energy difference varies in $x$, fronts may have to overcome a barrier, pointwise in $x$, in order to reduce energy in the system. For $\varepsilon>0$, one typically finds a \emph{pinning region} $(a_-,a_+)$ where two stationary interfaces exist, one stable and one unstable. At the boundary of this interval, the two stationary interfaces disappear in a saddle-node bifurcation. For values $a=a_++\mu$, $\mu>0$, the speed of the interface scales as 
\begin{equation}\label{e:1/2}
c\sim \mu^{1/2}.
\end{equation}
This can be formally (and more rigorously) understood as induced by the time that a moving interface spends near a saddle-node bifurcation, where time scales with $\mu^{-1/2}$. A more detailed analysis even yields prefactors in these expansions; see the discussion of lattice systems, below. 

% 
% pinning region and how it vanishes at eps=0, somewhwere else? cusp where two folds merge?

\paragraph{Pinning in discrete systems.}

Spatially discrete media can be viewed as an extreme limiting case of spatially periodic media; see for instance \cite{schvv}. They exhibit phenomena that are very similar. Consider therefore
\begin{equation}\label{e:lat}
\frac{\rmd}{\rmd t} u_i=\frac{d}{2}\left(u_{i+1}-2u_i+u_{i-1}\right)+f(u_i),\qquad i\in\Z.
\end{equation}
Stationary interfaces solve the two-term recursion
\[0=\frac{d}{2}\left(u_{i+1}-2u_{i}+u_{i-1}\right)+f(u_{i}).\]
Roughly speaking, the phenomena mirror the case of spatially periodic media: interfaces are stationary in the pinning region $a\in(a_-,a_+)$ and propagate with speed $c\sim \mu^{1/2}$ for $a-a_+=\mu>0$; see for instance \cite{carpio} for asymptotic and numerical studies,  \cite{mprev} for an overview of the theory, and \cite{mp1} for extensive and general theoretical results. 

% %They also concluded that a pinned front can begin to %propogate in two different manners. One such "depinning %transition" occurs when a continuous front profile becomes %discontinous as it approaches the pinning region. 
% 
% \paragraph{Continuous, periodic media}
% The pinned $u$ from the discrete case can be embedded as a function in $\mathbb{R}$ by using the equation:
%  \begin{equation*}
% 0=d\left(u(x+1)-2u(x)+u(x-1)\right)+f_a(u(x)).
% \end{equation*}
% The same pinned front exists, though now in a continuous medium. This form represents the kernel $\mathcal{K}=\frac{\delta(x-1)+\delta(x+1)}{2}$. By changing the nonlinearity $f_a(u)$, we can alter the shape and even existence of a pinning region in $(a,d)$-space.

\paragraph{Transversality in spatial dynamics.}

Pinned interfaces can be viewed as heteroclinic orbits of the two-term recursion, which defines a (local) diffeomorphism of the plane,
\[
u_{i+1}=u_i+w_i,\qquad w_{i+1}=w_i-\frac{1}{d}f(u_i+w_i).
\]
Indeed, $u=0,1$ and $w=0$ define hyperbolic fixed points of this diffeomorphism. The corresponding stable and unstable manifolds intersect along orbits that yield stationary interfaces. Such an intersection of stable and unstable manifolds in diffeomorphisms is typically (beyond this particular case)  transverse, hence robust with respect to changes in the parameter $a$. 

In a similar fashion, stationary profiles in spatially periodic media solve a non-autonomous differential equation 
\[
u_{x}=v,\qquad v_x=-u(1-u)(u-a+\varepsilon \sin(x)),
\]
whose time evolution $\Psi_{2\pi,0}$ defines a diffeomorphism of the plane with hyperbolic fixed points $(1,0)$ and $(0,0)$. Intersections of stable and unstable manifolds are typically transverse since time-translation symmetry is broken. 

In both cases, the boundary of the pinning region is given by a parameter value where stable and unstable manifolds intersect non-transversely, typically with a quadratic tangency that reflects the generic saddle-node bifurcation alluded to earlier. 

Summarizing, the traditional view of pinning associates open pinning regions with the absence of a continuous translational symmetry in the system.

% 
% 
% 
% Transversality is a condition that brings about the "jumping" behavior that we will observe when investigating wave speed asymptotics near the pinning region. Wave speed trajectories travel along a smooth manifold, then "jump" through a fast manifold at a "fold" point, and return to a slow manifold.  

\paragraph{Unpinning and speed asymptotics.}

\begin{figure}[h!]
\centering
\subfigure[Pinned.]{
\includegraphics[width=0.3\textwidth]{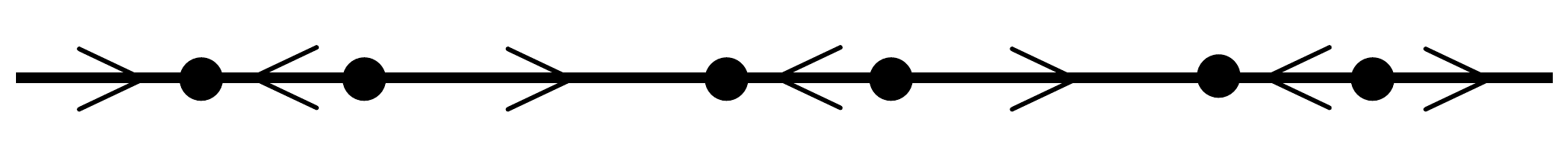}}
\subfigure[Boundary of pinning region.]{
\includegraphics[width=0.3\textwidth]{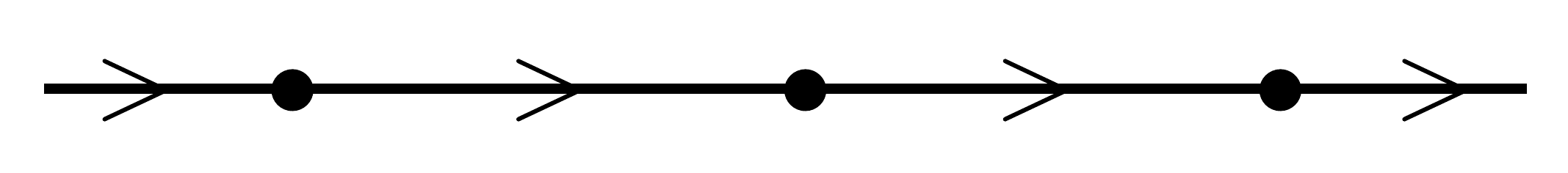}}
\subfigure[Unpinned.]{
\includegraphics[width=0.3\textwidth]{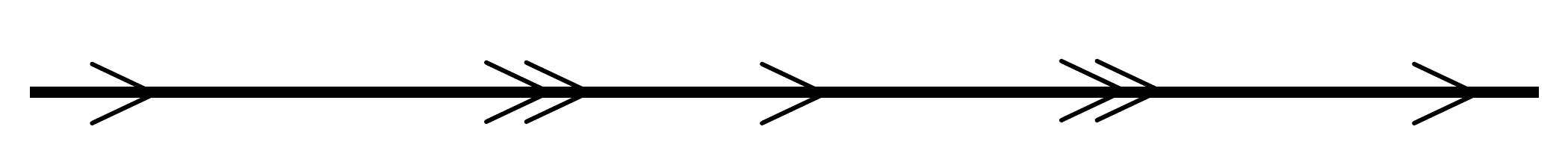}}
\caption{Temporal dynamics in function space, with the unfolding of a saddle-node and an associated homoclinic. Discrete spatial translation symmetry acts via shift to the right in this reduced phase space picture.}
\label{f:unpin}
\end{figure}

The previous discussion suggests that dynamics at the boundary of the pinning region are universal at leading order. The linearization at the stationary interface possesses a zero eigenvalue, associated with the saddle-node bifurcation (or the non-transverse intersection of stable and unstable manifolds, depending on the view point). The critical pinned front possesses a temporal homoclinic, corresponding to the translation of the interface by precisely one lattice site. Strictly speaking, this homoclinic orbit is homoclinic up to the discrete lattice translation symmetry. The unpinning transition then can be viewed as the unfolding of a saddle-node homoclinic orbit, with periodic orbits in the unpinned regime and heteroclinic orbits in the pinned regime; see Figure \ref{f:unpin}. 

\subsection{Discrete versus nonlocal systems}

We are interested here in, apparently quite different, nonlocal systems,
\begin{equation}\label{e:nonl}
u_t=d(-u+\mathcal{K}*u)+f(u),
\end{equation}
with a strongly localized convolution kernel $\mathcal{K}$. Here, $(\mathcal{K} \ast u)(x)=\int_{-\infty}^{\infty} K(x-y)u(y)\mathrm{d}y$ denotes  convolution. Existence and stability of interfaces in such systems has been studied extensively, under various assumptions on the kernel $\mathcal{K}$ and the nonlinearity $f$; see \cite{bates,chen}. As we shall explain in more detail later, it was noted that for weak coupling strength, $d\ll 1$, interfaces are discontinuous and do not propagate, for values of $a$ in a pinning region $(a_-,a_+)$. In other words, the presence of an open pinning region, here, is associated with a lack of regularity in the profile, rather than the absence of a translational symmetry. 

One can however emphasize similarities with lattice systems by embedding the lattice system (\ref{e:lat}) into a system on the real line,
\begin{equation}\label{e:disnonl}
\frac{\rmd}{\rmd t}u(x)=d\left(-u(x)+\frac{1}{2}(u(x+1)+u(x-1))\right)+f(u(x)).
\end{equation}
Of course, this system decouples in an infinite family of lattice systems $x\in x_*+\Z$, $x_*\in[0,1)$, each of which is equivalent to (\ref{e:lat}). While quite artificial, (\ref{e:disnonl}) exhibits the similarities between the different nonlocal and discrete pinning when written in the form (\ref{e:nonl}) with $\mathcal{K}(x)=\frac{1}{2}(\delta(x-1)+\delta(x+1))$ (although such kernels are not covered by assumptions in the references cited above). More explicitly, (\ref{e:disnonl}) possesses a continuous translational symmetry, but stationary interfaces are discontinuous, given for instance as $u(x)=u_{[x]}$, where $[x]$ is the integer part of $x$ and $u_j$ is the stationary interface in (\ref{e:lat})\footnote{In this sense, the profiles have countably many discontinuities at locations $x\in\Z$.}.

The point of view taken here is that (\ref{e:disnonl}) is a special element of the class of equations  (\ref{e:nonl}), in the sense that its kernel possesses very low regularity. One can then consider smoothed out versions of $\frac{1}{2}(\delta(x-1)+\delta(x+1))$ and ask about pinning regions and unpinning asymptotics. Our results indicate that both depend in a crucial fashion on the regularity of the approximation. While pinning is generic for the discrete kernel, pinning occurs only for sufficiently strong coupling in smooth kernels. Unpinning asymptotics are changed from speeds scaling with a exponent $1/2$ (\ref{e:1/2})  power law to 
\[
c\sim \mu^{3/2},
\]
in the pinning regime, or smooth speed asymptotics (\ref{e:1/1}) in the unpinned regime. We also find intermediate power laws 
\[
c\sim \mu^{5/4}.
\]
at critical coupling strengths.

\subsection{Summary of main results}

We consider nonlocal evolution equations of the form
\begin{equation}
 u_t=d(-u+\mathcal{K} \ast u) +f_a(u),\qquad x\in\R,\ u(t,x)\in\R.
\label{eq:a}
\end{equation}
Here, $f_a(u)$ is a bistable nonlinearity in $u$ with parameter $a$ detuning the energy levels of stable equilibria, and $d>0$ denotes coupling strength. We will mostly focus on  $f_a(u)=u(1-u)(u-a)$ but also comment on other cases.
We study traveling or stationary waves, setting $\xi:=x-c t$, and obtain
\begin{equation}
-c u_\xi=d\left(-u+\mathcal{K}\ast u\right)+f_a(u).
\label{eq:c}
\end{equation}

Our main results describe pinning regions and unpinning asymptotics as follows.
\paragraph{Pinning regions.} 
Open pinning intervals occur for $d<d_*$, $d_*>0$. In $(a,d)$-space, the pinning region $a\in (a_-(d),a_+(d))$ forms a cuspoidal region, 
\begin{equation}\label{e:cusp}
a_\pm(d)\sim a_*+a_2 (d-d_*)^2,
\end{equation}
with an explicit expression for cubic or sawtooth nonlinearities. 

\paragraph{Unpinning asymptotics.}
For $\mathcal{K}(x)=\frac{1}{2}\rme^{-|x|}$ and  $\mathcal{K}(x)=\rme^{-x}\chi_{[0,\infty)}(x)$, we establish rigorous asymptotics 
\begin{equation}\label{e:casy}
c=k_1|a-a_\pm|^{3/2}(1+\rmo(1)),
\end{equation}
and determine $k_1$ explicitly. We also outline a construction for asymptotics when $d=d_*$, and find
\begin{equation}\label{e:casyc}
c=k_1|a-a_\pm|^{5/4}(1+\rmo(1)),
\end{equation}
again with explicit $k_1$.

\paragraph{General kernels.} 
We study pinning asymptotics numerically, using \textsc{auto07p} and direct simulations. In particular, we corroborate our results numerically, confirming the value of $k_1$. We also show that our results carry over to a much larger class of smooth kernels, and we explore the dependence of scaling exponents in unpinning asymptotics (\ref{e:casy}) on smoothness of the kernel in a family of  kernels $\mathcal{K}$ with Fourier transform $\hat{\mathcal{K}}(\ell)=(1+\ell^2)^{-\beta/2}$.

\paragraph{Technical contribution.} 
We interpret the unpinning transition as a slow passage through a fold in a singularly perturbed system. For the special kernels mentioned above, the Fourier transform is rational and the nonlocal equation can be written as a system of ordinary differential equations. For small speeds, this system possesses a fast-slow structure, which we elucidate using geometric blowup methods. As a main result, we are able to give precise leading order expansions. We also explore the slow passage through an inflection point that, to our knowledge, has not been studied before.

\paragraph{Outline.} We state and derive our results on the shape of pinning regions in Section \ref{sec:pinning}. Section \ref{sec:asy} contains our main analysis in the case of rational kernels. Section \ref{sec:numerics} contains mostly numerical results, confirming the asymptotics (in a surprisingly large parameter regime), and exploring other kernels and scalings. We conclude with a brief discussion in Section \ref{sec:dis}.

\section{Pinning regions}\label{sec:pinning}

The phenomenon of pinning was noticed and precisely characterized in \cite{bates}; see also \cite{chen} for more general statements. In this section, we state the result from \cite{bates} formally and characterize pinning regions in several explicit cases and in a generic example. 

Roughly speaking, when characterizing pinned fronts, one notices that, for positive kernels, the system possesses a comparison principle and finds that interfaces are monotone. Stationary interfaces solve 
\[
0=d(-u+\mathcal{K}\ast u)+f_a(u).
\]
Since $\mathcal{K}*u$ is monotone and smooth, provided the kernel is sufficiently smooth, the remainder 
\begin{equation}
g(u;a,d):=u-\frac{1}{d}f_a(u)
\label{eq:b}
\end{equation}
is monotone and smooth. This is however not possible when $g$ is not monotone. On the other hand, non-smooth profiles cannot propagate due to the continuity of the evolution of (\ref{eq:a}). Moreover, profiles depend continuously on $a$ in $L^\infty$ (under an open set of conditions), so that discontinuous, pinned profiles exist for open pinning regions.

\begin{figure}[h!]
\centering
\begin{minipage}{.5\textwidth}
\centering
\includegraphics[width=.6\linewidth]{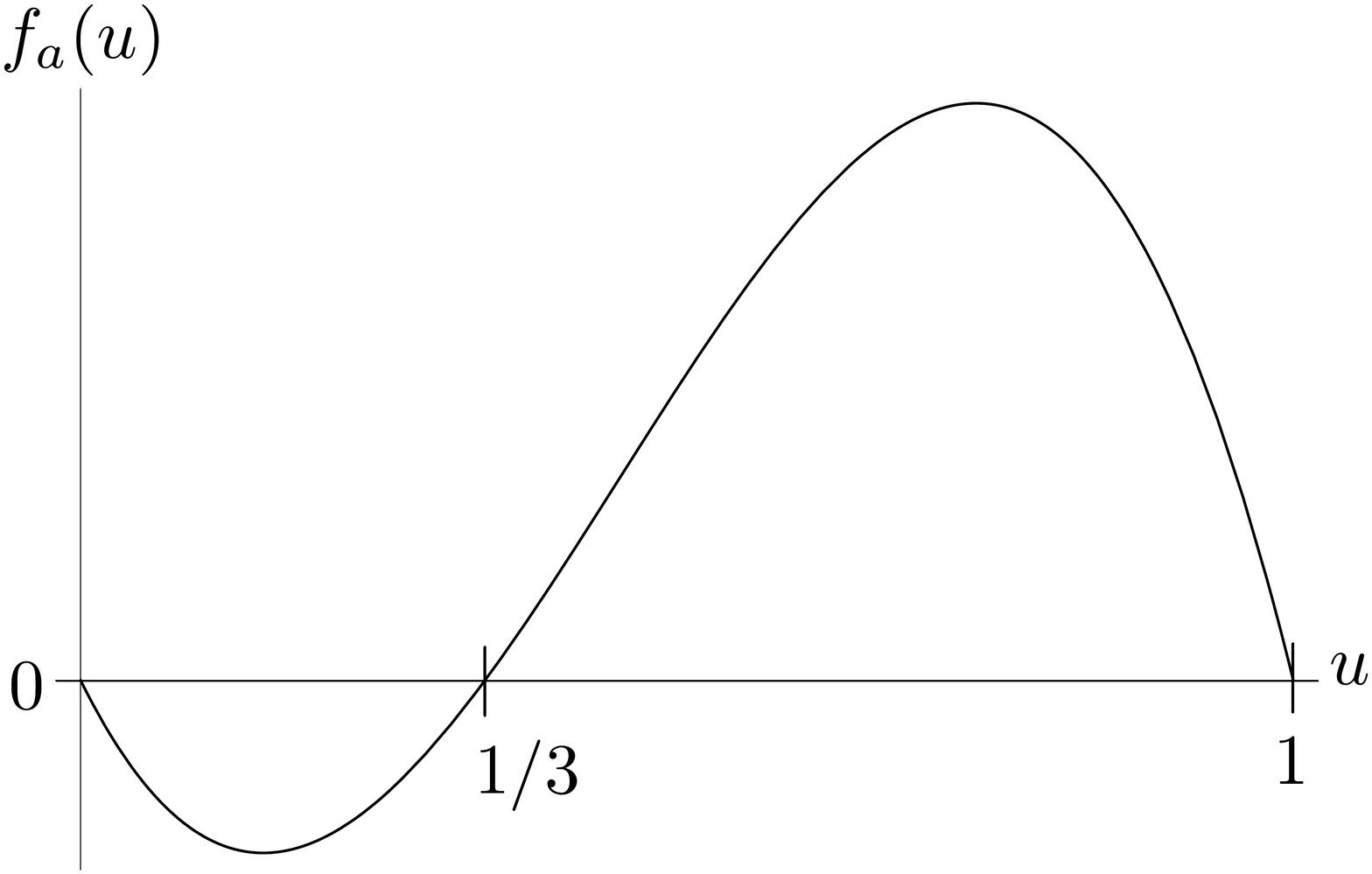}
\end{minipage}%
\begin{minipage}{.5\textwidth}
\centering
\includegraphics[width=.6\linewidth]{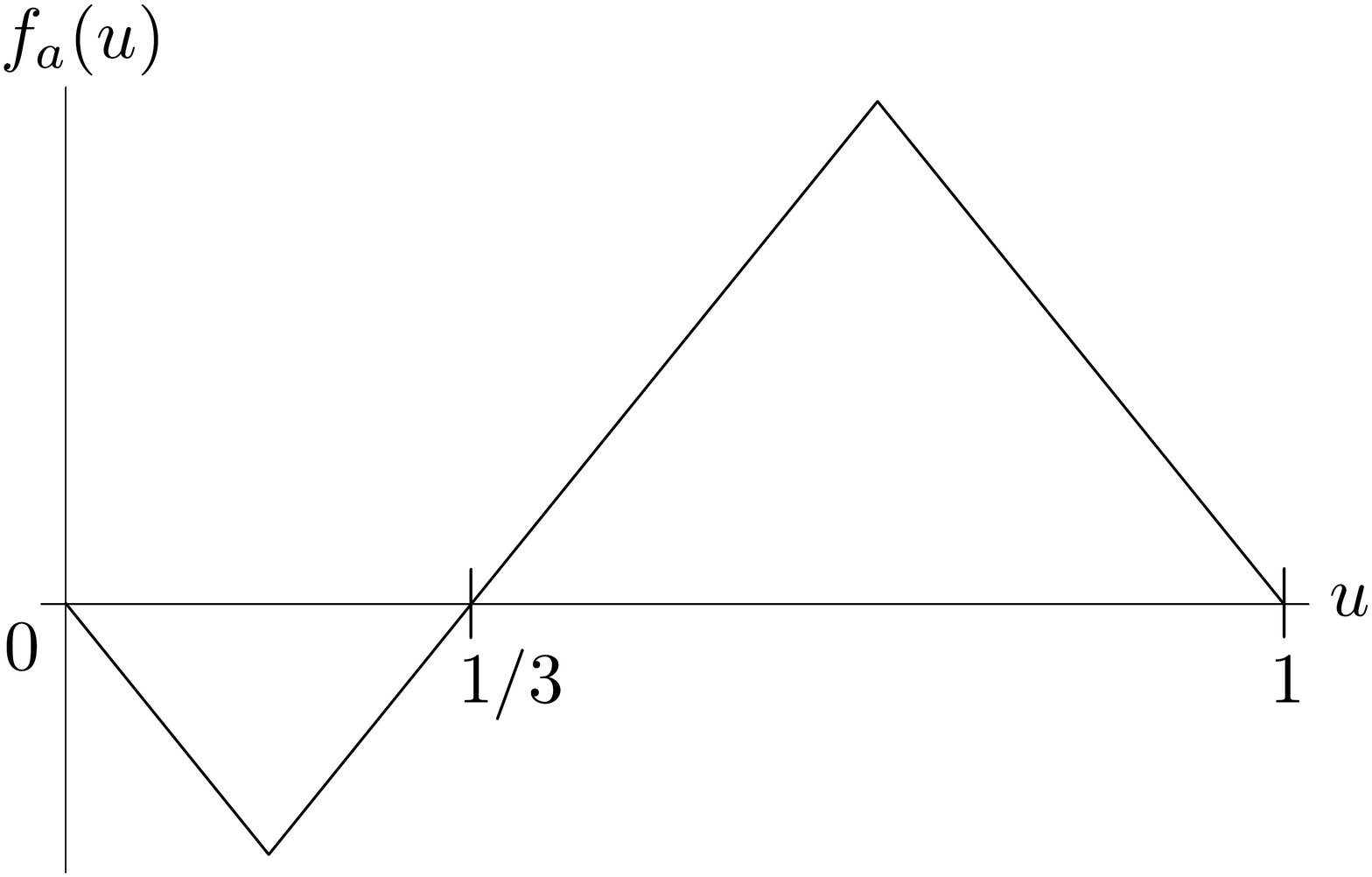}
\end{minipage}
\caption{Cubic (left) and piecewise linear (right) nonlinearity with $a=\frac{1}{3}$.}
\label{fig:nonlinearities}
\end{figure}

To be more precise, we revisit the setting of \cite{bates} in more detail.
\begin{Hypothesis}{(H1)} We require the following conditions for the nonlinearity $f_a(u)$:
\begin{enumerate}
\item $f_a(u) \in C^3(\mathbb{R})$;
\item $f(u)=0$ precisely when $u\in\{0,a,1\}$, where $0<a<1$;
\item $f'_a(0),f'_a(1)<0$.
 \end{enumerate}
 \label{hyp:h1}
\end{Hypothesis}
Figure \ref{fig:nonlinearities} offers two examples of functions $f_a(u)$ for which we will calculate the pinning regions later in this section. We note that the cubic nonlinearity satisfies (H1). The piecewise linear nonlinearity on the right does not satisfy all criteria listed in (H1), but we shall explain later that the results from \cite{bates} still apply to this function.
% 
% \begin{figure}[h!]
% \centering
% \subfigure[]{
% \includegraphics[width=0.35\textwidth]{Figures/cubiconethird}}
% \hspace{.4in}
% \subfigure[]{
% \includegraphics[width=0.35\textwidth]{Figures/piecewiselinearonethird}}
% \caption{Cubic (left) and piecewise linear (b) nonlinearity with $a=\frac{1}{3}$.}
% \label{fig:nonlinearities}
% \end{figure}

\begin{Hypothesis}{(H2)}
We assume that $g(u;a,d)$ from \eqref{eq:b} is monotone on a maximum of 3 intervals.
\end{Hypothesis}
In other words, 
\begin{equation*}
g'(u;a,d)>0 \text{ for } u \in [0,\beta_1) \cup (\beta_2,1], \quad  g'(u;a,d)<0 \text{ for } u \in (\beta_1,\beta_2),
\end{equation*}
for $0<\beta_1 \leq \beta_2<1$. When $\beta_1=\beta_2$, $g(u;a,d)$ is monotone on $[0,1]$ and $g'(\beta_1;a,d)\geq 0$.  We can now define
\begin{eqnarray*}
 g_k(u;a,d)=\left\{ \begin{array} {ll}
g(u;a,d) & \mbox{$u\in[0,\beta_1]\cup[\beta_2,1]$,} \\
k & \mbox{$u \in [\beta_1,\beta_2]$},  \\
\end{array}
\right.
\end{eqnarray*}
where $g(\beta_1;a,d)=g(\beta_2;a,d)=k$.

\begin{Hypothesis}{(H3)}
We require that the convolution kernel $\mathcal{K}$ satisfies:
\begin{enumerate}
\item $\mathcal{K}, \mathcal{K}' \in L^1(\mathbb{R})$ with $\int \mathcal{K}(x)\rmd x= 1$;
\item $\mathcal{K}(x)=\mathcal{K}(-x) \geq 0$;
\item $\int \mathcal{K}(x)|x|\rmd x < \infty$.
\end{enumerate}
\end{Hypothesis}
We are now ready to recall the main result that our analysis of pinning regions is based on.
\begin{theorem}\cite{bates}
Assume (H1)-(H3). Let  $(c,u)$ be a solution of equation \eqref{eq:c}, then:
\begin{enumerate}
\item $u$ has at most one jump discontinuity;
\item the solution $(\eps,u)$ is unique up to translation in the class of monotone profiles $u$;
\item if $c \neq 0$, then $u\in \mathcal{C}^4$;
\item  $c=0$ if and only if  $\int_0^1 g_k(v;a,d)\mathrm{d}v =\frac{1}{2}$ for some $k$.
\end{enumerate}
\label{thm:bates}
\end{theorem}

In the remainder of this section, we shall analyze four examples, two of which are covered by Theorem \ref{thm:bates}.

\subsection{Cubic nonlinearity}\label{subsec:cubicpin}

With these criteria in mind we proceed by calculating the pinning region in the simple cubic case  $f_a(u):=u(1-u)(u-a)$, applying Theorem \ref{thm:bates} in a straight forward fashion.

\begin{figure}[h!]
\centering
\subfigure[]{
\includegraphics[width=0.45\textwidth]{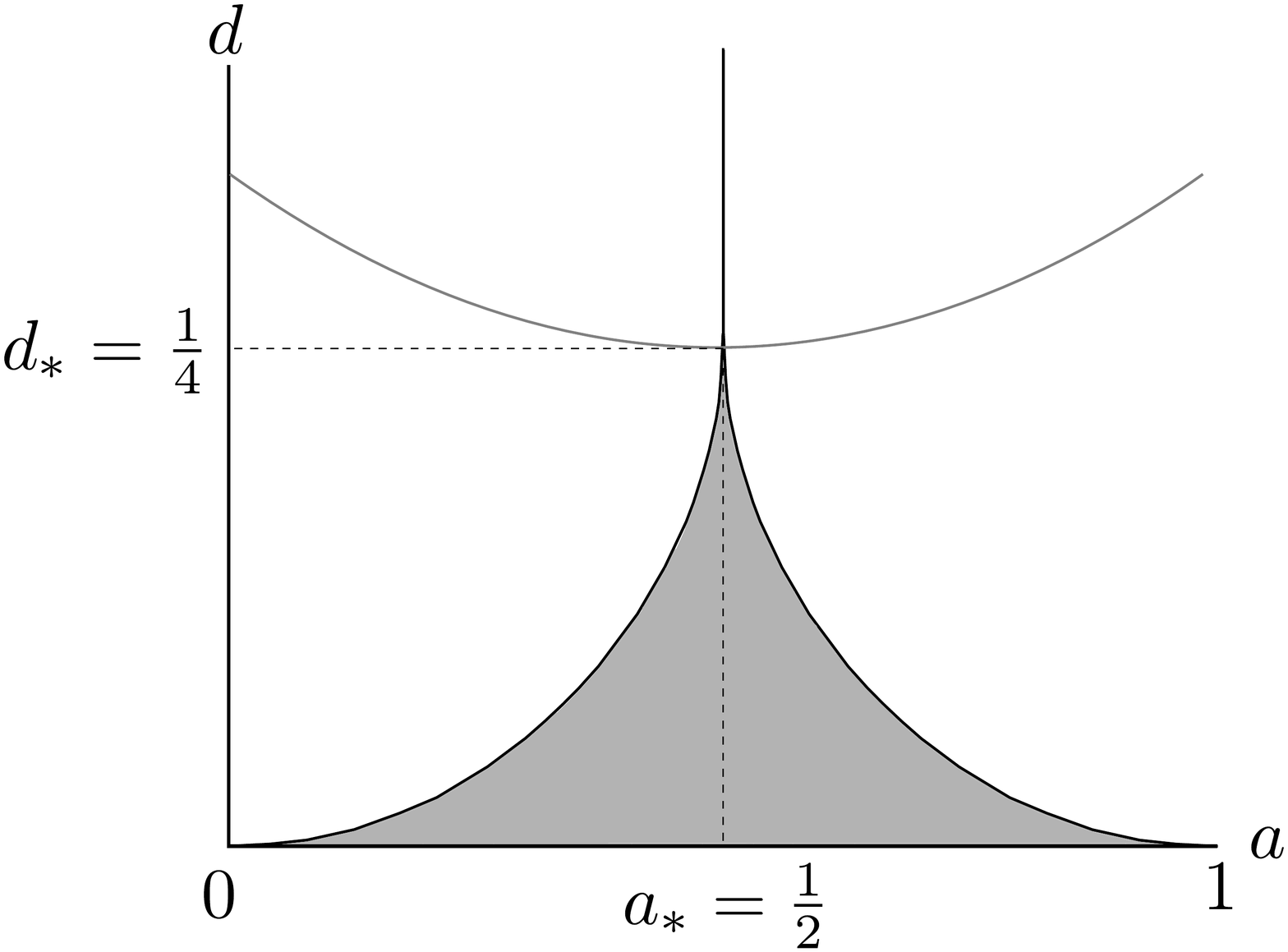}}
\hspace{.1in}
\subfigure[]{
\includegraphics[width=0.45\textwidth]{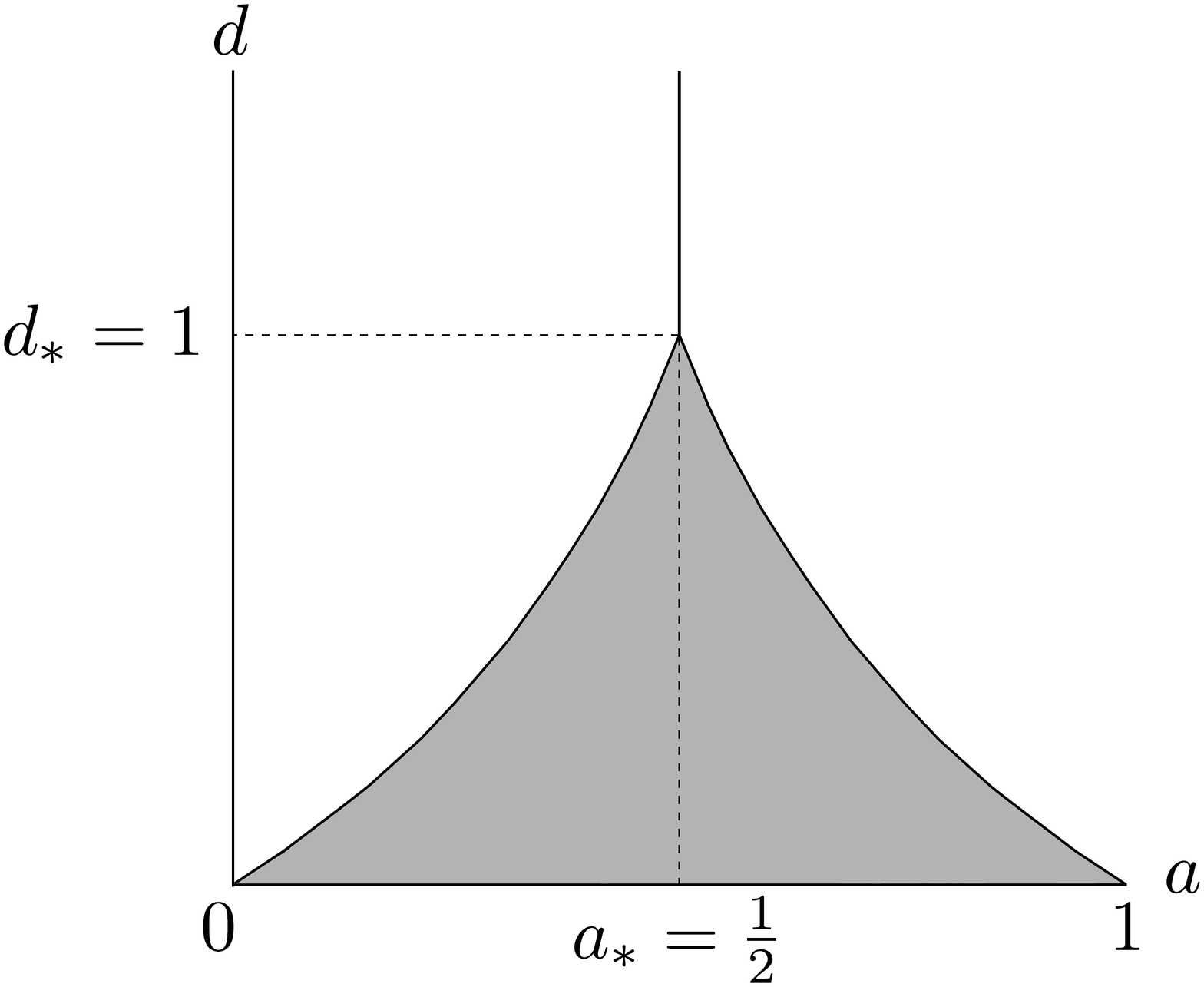}}
\caption{Pinning regions for cubic (a) and piecewise linear (b) functions. }
\label{fig:generalcubicpinningregion}
\end{figure}

\begin{lemma}
The pinning region of the cubic  $f_a(u)=u(1-u)(u-a)$ in the $(a,d)$-plane is bounded by
\begin{eqnarray*}
 d(a)=\left\{ \begin{array} {ll}
\frac{1}{3}(1-a+a^2-\sqrt{1-2a})& \mbox{$a\leq \frac{1}{2}$,} \\
\frac{1}{3}(1-a+a^2-\sqrt{-1+2a})& \mbox{$a \geq \frac{1}{2}$.}  \\
\end{array} 
\right.
\end{eqnarray*}
In particular, robust pinning occurs for $d<1/4$, in an interval $(a_-(d),a_+(d))$, with 
\begin{equation}
a_\pm(d)=\frac{1}{2}\pm \frac{9}{2}\left(d-\frac{1}{4}\right)^2+\rmO\left(\left(d-\frac{1}{4}\right)^4\right);
\label{e:apm}
\end{equation}
see Figure \ref{fig:generalcubicpinningregion}.
\end{lemma}

\begin{Proof} We construct $g_k$ as defined after (H2); see Figure \ref{fig:generalcubic}. Let  $u_-^M$ denote the location of the local maximum of $g(u;a,d)$ and $u_+^M$ the other preimage of this local maximum. Similar, define $u_-^m$ and $u_+^m$ using the local minimum. Now define
\begin{eqnarray}
 g_M(u;a,d):=\left\{ \begin{array} {ll}
g(u;a,d) & [0,u_-^M]\cup [u_+^M,1], \\
M&(u_-^M, u_+^M),
\end{array} 
\right.
\label{eq:gmod}
\end{eqnarray}
and $g_m$ in an equivalent fashion.

\begin{figure}[h!]
\centering
\begin{minipage}{.5\textwidth}
 \centering
\includegraphics[width=105mm]{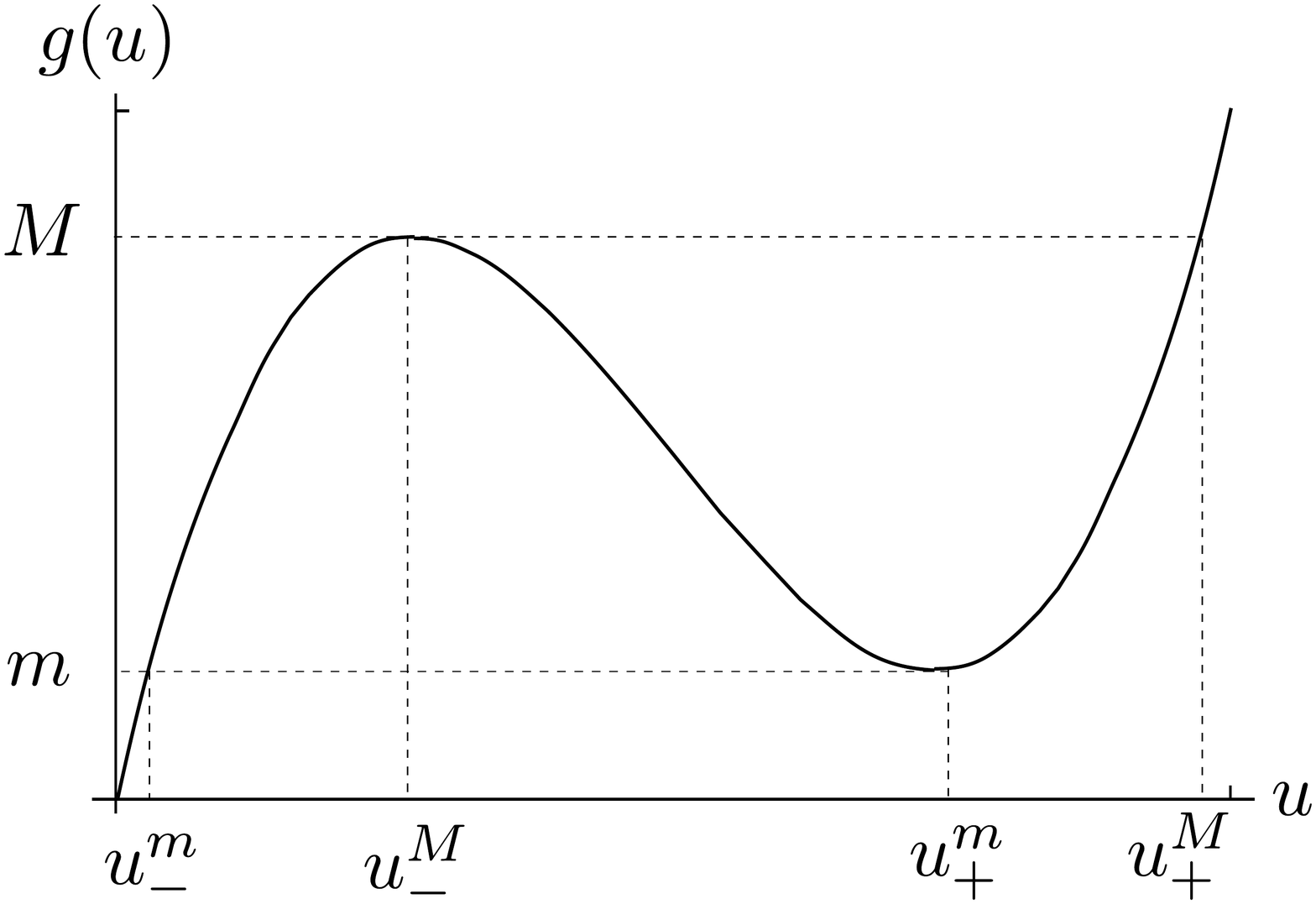}
\end{minipage}%
\begin{minipage}{.5\textwidth}
  \centering
  \includegraphics[width=.6\linewidth]{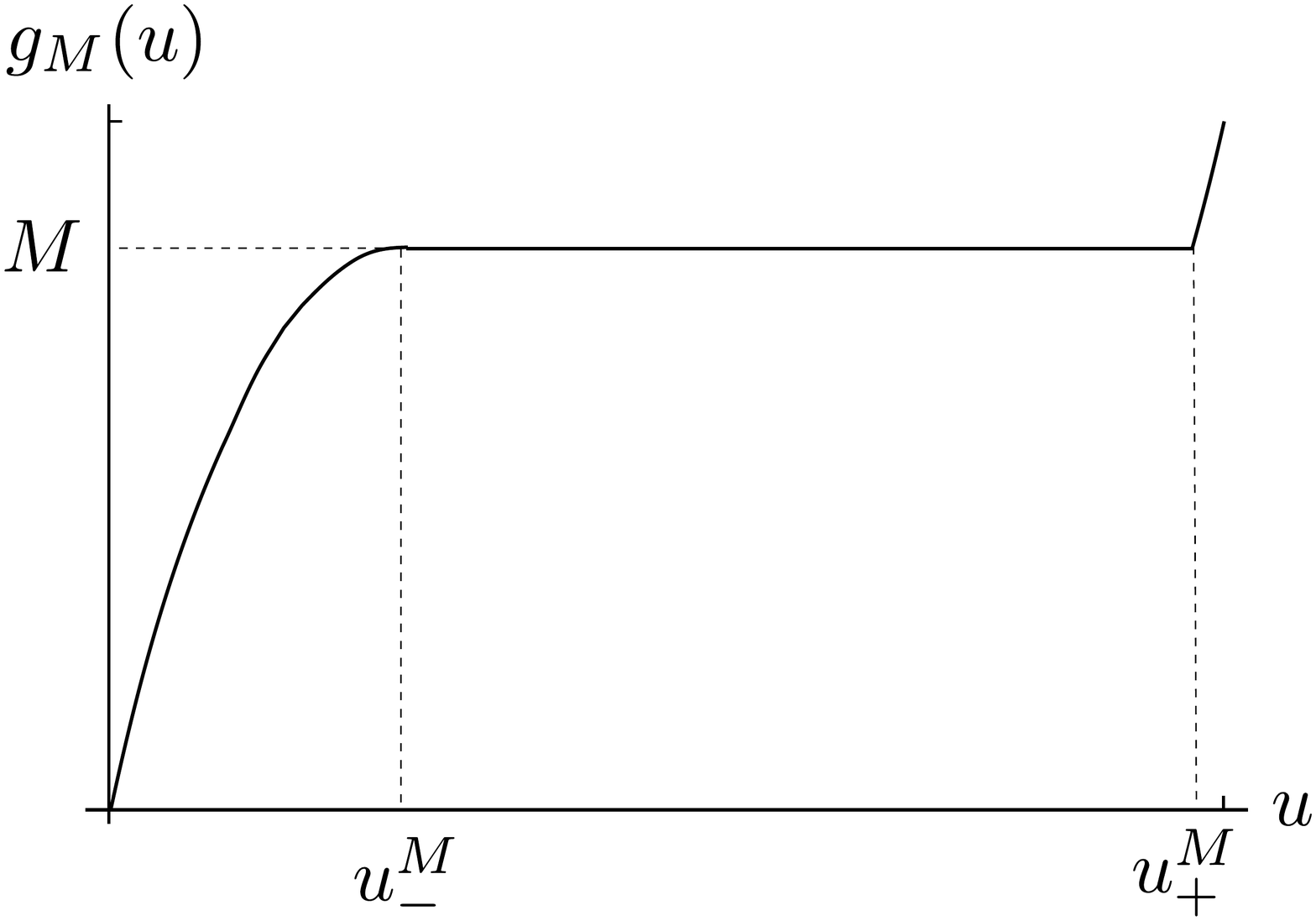}
 % \label{fig:test1}  
  \centering
  \includegraphics[width=.6\linewidth]{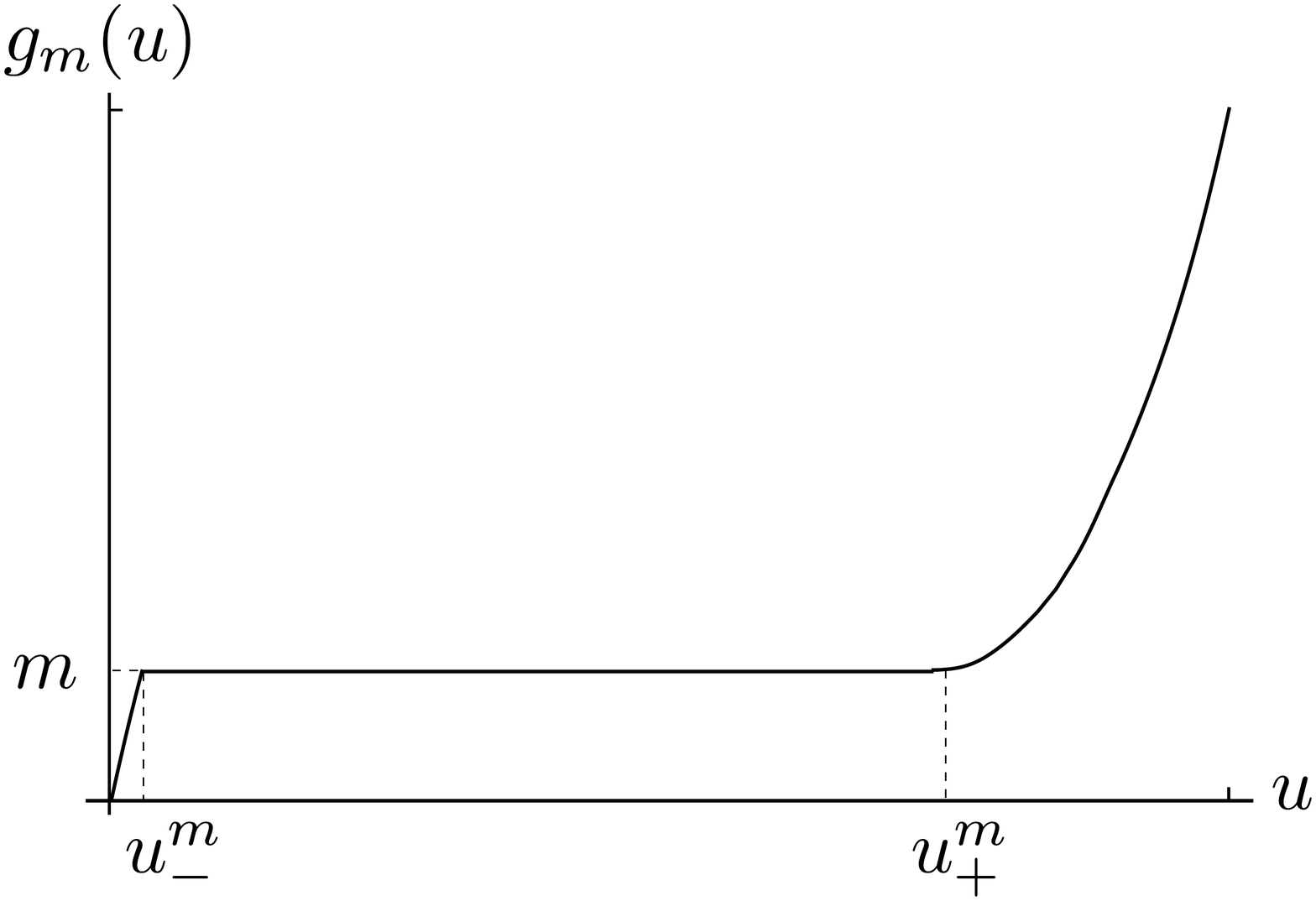}
%  \label{fig:test2}
\end{minipage}
\caption{Cubic nonlinearity $g(u)$ and modified functions $g_M(u)$  and $g_m(u)$; see \eqref{eq:gmod}.}
\label{fig:generalcubic}
\end{figure}

From Theorem \ref{thm:bates}, (iv), we  find that the conditions  $\int_0^1 g_{m/M}(v;a,d)\mathrm{d}v=1/2$ determine the boundaries of the pinning region. We can explicitly compute those integrals and solve for $d$, to obtain 
\begin{equation}
d_+=\frac{1}{3}\left(1-a+a^2-\sqrt{1-2a}\right),\qquad 
d_-=\frac{1}{3}\left(1-a+a^2-\sqrt{-1+2a}\right).
\label{eq:cub_reg}
\end{equation}
\end{Proof}

%\begin{figure}[h!]
%\centering
%\includegraphics[height=2.5in]{Figures/GeneralCubicPinningRegion}\hspace*{10mm}
%\includegraphics[height=2.5in]{Figures/PinningRegionSeesaw}
%\caption{Pinning regions for cubic (left) and piecewise linear (right) functions. }
%\label{fig:generalcubicpinningregion}
%\end{figure}

\subsection{Generic nonlinearity}\label{subsec:genpin}

We show that the shape of the pinning region is universal, scaling $a_\pm(d)\sim (d-d_*)^2$. We therefore make assumptions on a nonlinearity in ``general position'' near the tip of the pinning region. 

\begin{Hypothesis}{(H4)} We require that $g(u;a,d)$ is smooth and satisfies for $(a,d)\sim (a_*,d_*)$,
\begin{enumerate}
	\item \emph{area balance}: $\int_0^1 g(v;a_\ast,d)\rmd v=\frac{1}{2}$;
	\item \emph{inflection point}: $g'(u_\ast;a_\ast,d_\ast) = g''(u_\ast;a_\ast,d_\ast) = 0$, and $g'''(u;a,d)>0$;
	\item\emph{$d$-unfolding}: $\partial_{d,u}g(u_\ast;a_\ast,d_\ast)>0$;
	\item \emph{$a$-unfolding}: $\partial_a\int_0^1  g(v;a_\ast,d_\ast)\mathrm{d}v\neq 0$.
\end{enumerate}
For convenience, in accordance with the specific form of $g$ in (\ref{eq:b}), we also assume that $\partial_d g(u_\ast;a_\ast,d)=0$ for all $d\sim d_\ast$.
\end{Hypothesis}
Before stating the main result of this section, we define the shorthand notation $g_3:=\partial_u^3g(u_\ast;a_\ast,d_\ast)$, $g_{u,d}:=\partial_d\partial_u g(u_\ast;a_\ast,d_\ast)$ and $g_\ast:=g(u_\ast;a_\ast,d_\ast)$.

\begin{theorem}
Assuming (H1)-(H4), the pinning region in a neighborhood of $(a_*,d_*)$ is contained in $d<d_*$, of the form $(a_-(d),a_+(d))$, with smooth functions $a_\pm(d)$ and expansions
\begin{equation*}
\begin{split}
a_\pm(d)&=a_\ast\pm a_2(d_\ast-d)^2+\mathcal{O}\left((d_\ast-d)^3\right),\\
a_2 &=-\frac{9}{2}\frac{g_{u,d}}{g_3\int_0^1\partial_a g(v;a_\ast,d_\ast)\mathrm{d}v}.\\
\end{split}
\end{equation*}
\label{thm:pinningregion}
\end{theorem}
\begin{Proof}
We focus on finding the lower boundary of the pinning interval, associated with the local maximum of $g$ at $u=u_-^M$.  To simplify our notation, we drop the superscript $M$ and write $u_-$ for the location of the maximum and $u_+$ for the second preimage of the maximum value. We also translate our system to the critical $u$, $d$, and $a$ values: $\tilde{u}=u-u_\ast$, $\lambda=d_\ast-d$ and $\mu=a_\ast-a$.

We  expand $g$ in $u$, $a$ and $d$ around $u_\ast$, $a_\ast$ and $d_\ast$. We introduce the scaling $\lambda=\eps^2$, $\mu=\eps^4\mu_1$, and $\tilde{u}=\eps u_1$.
%and scale $u_\ast$ with $\eps$ as well. Then,
Then,
\begin{equation*}
g(u_\ast+\eps u_1;a_\ast-\eps^4\mu,d_\ast-\eps^2)= g_\ast +\frac{1}{6} g_3\eps^3 u_1^3-g_{u,d}\eps^3 u_1+\mathcal{O}\left( \eps^4\right).
\end{equation*}
In this scaling, we find the location of the local maximum, bifurcating from the inflection point,  from 
\[
0=g_u(u_\ast+\eps u_1;a_\ast-\eps^4\mu,d_\ast-\eps^2)= \eps^3\left( \frac{1}{2}g_3 u_1^2-g_{u,d}\right) + \mathcal{O}\left(\eps^3\right),
\]
which gives the maximum location $u_{1,-}(\eps,\mu)$ as
\[u_{1,-}=-\sqrt{\frac{2g_{u,d}}{g_3}}+\mathcal{O}(\eps).\]
The value $g_-(\eps,\mu)$ at the maximum is 
\begin{equation*}
% \begin{split}
g(u_\ast+\eps u_{1,-};a_\ast-\eps^4\mu,d_\ast-\eps^2)= g_\ast+\eps^3\left(\frac{1}{6}g_3 u_{1,-}^3 -g_{u,d}u_{1,-}\right)+\mathcal{O}(\eps^4)
=g_\ast + \eps^3\sqrt{\frac{8g_{u,d}^3}{9g_3}}+\mathcal{O}(\eps^4).
% \end{split}
\end{equation*}

Finally, we find  $u_{+,1}(\eps,\mu)$, solving 
\[
g(u_\ast+\eps u_{1,-};a_\ast-\eps^4\mu,d_\ast-\eps^2)=
g(u_\ast+\eps u_{1,+};a_\ast-\eps^4\mu,d_\ast-\eps^2),
\]
which readily gives $u_{1,+}=-2u_{1,-}+\rmO(\eps)=\sqrt{\frac{8g_{u,d}}{g_3}}+\mathcal{O}(\eps)$.

Now, we can evaluate the integral of  $g_M(u;a,d)$ near $a_\ast$ and $d_\ast$:
\begin{equation*}
\begin{split}
\int_0^1 g_M(v;a,d)\mathrm{d}v&=\int_0^1 g(v;a,d)\mathrm{d}v-\int_{u_\ast+\eps u_{1,-}} ^{u_\ast+\eps u_{1,+}}[g(v;a,d)-g(u_{1,-};a,d)]\mathrm{d}v,\\
&=\frac{1}{2}+\eps^4 \mu_1\int_0^1 \partial_a g(v;a_\ast,d_\ast)\rmd v -\eps  \int_{u_{1,-}}^{u_{1,+}}[g(u_*+\eps v_1;a_*,d_*-\eps^2)-g_-(\eps,\mu)]\rmd v_1+\mathcal{O}\left(\eps^5\right).\\
\end{split}
\end{equation*}
The second integral can be further evaluated as
\begin{equation*}
\begin{split}
\int_{u_{1,-}}^{u_{1,+}}[g(u_*+\eps v_1;a_*,&d_*-\eps^2)-g_-(\eps,\mu)]\rmd v_1\\
&=\eps^3\int_{u_{1,-}}^{u_{1,+}} \left[\frac{1}{6} g_3v_1^3-g_{u,d} v_1 - \sqrt{\frac{8g_{u,d}^3}{9g_3}}\right]\mathrm{d}v_1 + \mathcal{O}\left(\eps^4\right),\\
&=\eps^3\left(\frac{g_3}{24}(u_{1,+}^4-u_{1,-}^4) -\frac{g_{u,d}}{2}(u_{1,+}^2-u_{1,-}^2)-\frac{2}{3} \sqrt{\frac{2g_{u,d}}{g_3}}(u_{1,+} - u_{1,-}) \right)+\mathcal{O}\left(\eps^4\right),\\
&=\eps^3\left(\frac{15}{24}-3-4\right)\frac{g_{u,d}^2}{g_3}+\mathcal{O}(\eps^4)=-\eps^3 \frac{9g_{u,d}^2}{2g_3}+\mathcal{O}(\eps^4).\\
\end{split}
\end{equation*}
where in the third equality we substituted $u_{1,-} = -\sqrt{\frac{2g_{u,d}}{g_3}}+\mathcal{O}(\eps)=-\frac{u_{1,+}}{2}+\mathcal{O}(\eps)$. The pinning boundary is obtained when $\int g_M=1/2$, which gives
\begin{equation}
-\frac{9}{2}\frac{g_{u,d}^2}{g_3}\eps^4+\mathcal{O}(\eps^5)=\eps^4\mu_1\int_0^1 \partial_a g(v;a_\ast,d_\ast)\mathrm{d}v,
\end{equation}
or, in unscaled variables,
\begin{equation}
-\frac{9}{2}\frac{g_{u,d}^2}{g_3}\lambda^2+\mathcal{O}\left(\lambda^\frac{5}{2}\right) = \mu\int_0^1 \partial_a g(v;a_\ast,d)\mathrm{d}v.
\label{eq:cusp}
\end{equation}
Using opposite scalings, such as $\mu=-\eps^4$, $\lambda=\lambda_1\eps^2$, one finds that the pinning boundary is unique in $\mu<0$, which confirms the scaling and concludes the proof.
% 
% From this relationship we see that our original scaling choice in $\eps$ was correct. This calculation gives us a unique pinning boundary in $(a,d)$-space below $\lambda^2\sim\mu$. To show no pinning region exists above this curve simply repeat the above calculations with the scaling $\mu=\eps^4$, $\lambda=\lambda_1 \eps^2$ and see that the same solutions comes forth.
\end{Proof}

\subsection{Piecewise linear case}\label{subsec:plinpin}
Piecewise linear nonlinearities provide interesting explicit test cases. We consider here
\begin{eqnarray}
 f_a(u)=\left\{ \begin{array} {ll}
-u& \mbox{$u\leq \frac{a}{2}$,} \\
u-a& \mbox{$\frac{a}{2}\leq u \leq\frac{1+a}{2}$,}  \\
1-u& \mbox{$\frac{1+a}{2}\leq u$.}  \\
\end{array} 
\right.
\label{eq:seesaw}
\end{eqnarray}
Going through the proof in \cite{bates}, one can verify that the pinning criterion from Theorem \ref{thm:bates} (iv), is applicable although this function is not sufficiently smooth.
In the following, we simply apply this criterion formally.

The local maximum  is located at $u_-=\frac{a}{2}$ and  $f_a(\frac{1+a}{2})=f_a(u)$. Using the integral condition to solve for $d$ in terms of $a$, we find that the pinning region in the case of the piecewise linear function is bounded by:
\begin{eqnarray}
 d(a)=\left\{ \begin{array} {ll}
\frac{a}{1-a}& \mbox{$a\leq \frac{1}{2}$},\\
\frac{1-a}{a}& \mbox{$a\geq \frac{1}{2}$}.  \\
\end{array} 
\right.
\label{eq:e}
\end{eqnarray}
Note that the tip of the pinning region opens up at a finite angle rather than a cusp, in contrast to the smooth cubic and generic case; see Figure  \ref{fig:generalcubicpinningregion} for a comparison.

\subsection{Asymmetric, 1st order kernel : one-sided pinning}\label{subsec:1dpin}
Deviating from the setup in \cite{bates}, we consider the asymmetric kernel with Fourier symbol 
\begin{equation}
\label{eqKernel1d}
\hat{\mathcal{K}}(\ell)=\frac{1}{1+i\ell}.
\end{equation}
The equation for interfaces can then be written as an ordinary differential equation, 
\begin{equation}
\begin{split}
-c\,u_\xi&=d(w-u)+f_a(u),\\
w_\xi&=u-w.
\end{split}
\end{equation}
For $c=0$, the first equation is of  the form $w=g(u;a,d)$ with inverse $u=g_\pm(w;a,d)$ where $g_\pm$ refers to the smooth inverses, defined as continuations from $w\sim \infty$ and $w\sim -\infty$ respectively. We are looking for a solution where $w$ is continuous, and $u,w\to 1$ for $\xi\to-\infty$ and $u,w\to 0$ for $\xi\to\infty$.

Since $w=1$ is stable in $w_x=g^{-1}_+(w;a,d)-w$ and $w=0$ is stable in $w_x=g^{-1}_-(w;a,d)-w$, such solutions only exist when $1\leq g(u_{max};a,d)$. The boundary of the pinning region is therefore determined by the condition 
that $g(u;a,d)$ has a maximum at $1$, which gives 
$d=\frac{a^2}{4}.$
Similarly, the pinning condition for interfaces with  $u,w\to 1$ for $\xi\to\infty$ and $u,w\to 0$ for $\xi\to-\infty$ turns out to be $d=\frac{(1-a)^2}{4}.$ Summarizing, waves that connect 1 to 0 are pinned in 
\begin{equation*}
\left\{(a,d) \mid d\leq \frac{a^2}{4}, a\in (0,1)\right\},
\end{equation*}
while waves that connect 0 to 1 are pinned in 
\begin{equation*}
\left\{(a,d)\mid d\leq \frac{(1-a)^2}{4}, a\in (0,1)\right\}.
\end{equation*}
We will later construct traveling waves outside of this pinning region. 

\section{Speed asymptotics --- rational kernels}\label{sec:asy}
In this section, we derive asymptotic expansions for the speed of moving interfaces, close to but outside of the pinning region. Throughout, we fix $d$ and vary $a$ near $a_-$, the left boundary of the pinning region. We also exploit the fact that kernels where $\hat{\mathcal{K}}$ is rational, can be expressed by solving a linear differential equation. The section is split into three parts. We investigate the kernel $\hat{\mathcal{K}}=(1+\ell^2)^{-1}$, first in the generic case, and then at the tip of the pinning region; Sections \ref{subsec:speedsecondderivative} and \ref{subsec:prtip}. We then study $\hat{\mathcal{K}}=(1+\rmi\ell)^{-1}$ in Section \ref{subsec:onederivative}.

\subsection{Speed asymptotics for second derivative kernel}
\label{subsec:speedsecondderivative}
Consider 
$\mathcal{K}(x)=\rme^{-|x|}/2,$
with Fourier transform $\mathcal{\widehat{K}}(\ell)=\frac{1}{1+\ell^2}$. Writing $w=\mathcal{K}\ast u$, we see that
$ u=w-w''
$, so that  \eqref{eq:c} can be rewritten in the form
\begin{equation}
\begin{split}
w_x&=v,\\
v_x&=w-u,\\
-c u_x&=d(-u+w)+f_a(u).\\
\end{split}
\label{eq:sa_3ODE}
\end{equation}
For simplicity, we also focus on the cubic nonlinearity, which will allow us to explicitly determine constants for the leading order term in the expansion. 
\begin{theorem}
For fixed $d<d_\ast$, as the parameter $a$ approaches the left boundary of the pinning region at $a_-$ the asymptotics of the wave speed of the unique trajectory given in Theorem \ref{thm:bates} are:
\begin{equation*}
c=k_1(a_--a)^{\frac{3}{2}}+\mathcal{O}\left((a_--a)^2\ln(a_--a)\right), \text{ as } a\nearrow a_-.
\end{equation*}
The explicit formula for $k_1$ can be found in  \eqref{eq:k1}. Note that $a_-$  depends on $d$; see \eqref{e:apm}. The equivalent result (with same constant) holds for the right boundary of the pinning region.
\label{thm:2D}
\end{theorem}
\begin{Proof}
We want to analyze heteroclinic solutions that connect the saddle equilibria $\underline{u}_+=(0,0,0)$ and $\underline{u}_-=(1,0,1)$ in $\underline{u}=(w,v,u)$-space. With the parameter $c$, the system has a natural slow-fast structure. We will find that, at the pinning boundary, there exists a \emph{singular trajectory}, patched together from solutions of (\ref{eq:sa_3ODE}) at $c=0$, and a singular fast jump through a fold point. In order to identify this singular trajectory, we analyze (\ref{eq:sa_3ODE}) in slow and fast time, separately.

\paragraph{Fast system.} We first rescale time, $x=cy$, and find 
\begin{equation}
\begin{split}
w_y & = c v,\\
v_y & = c(w-u),\\
-u_y & = d(-u+w)+f_a(u).
\end{split}\label{e:fast}
\end{equation}
At $c=0$, this system possesses a manifold of equilibria $\mathcal{M}$ given through $w=g(u;a,d)$. We can write this manifold as a graph over the $(w,v)$-plane, up to certain fold points, by inverting the cubic $g$, writing $u=g_\pm^{-1}(w;a,d)$, where the inverses are understood as continuations from $\pm\infty$, respectively. In the following, we think of these manifolds as graphs with $u$ plotted in the vertical direction, so that $g_+^{-1}$ gives the upper branch $\mathcal{M}_+$ and $g_-^{-1}$ the lower branch $\mathcal{M}_-$ of the manifold; see Figure \ref{fig:cubicsurface} for an illustration of the slow manifold.  Except for  the fold points, both upper and lower branch of the  manifold of equilibria are normally hyperbolic in the terminology of \cite{fenichel:79}, therefore persist as locally invariant manifolds $\mathcal{M}=\mathcal{M}(c)$ for small $c\neq 0$. In fact, both are unstable in the normal direction\footnote{The middle branch is normally hyperbolic except for the fold points, as well, 
and normally stable, but of less relevance to us, here.}.
\begin{figure}[h]
\centering
\includegraphics[scale=0.35]{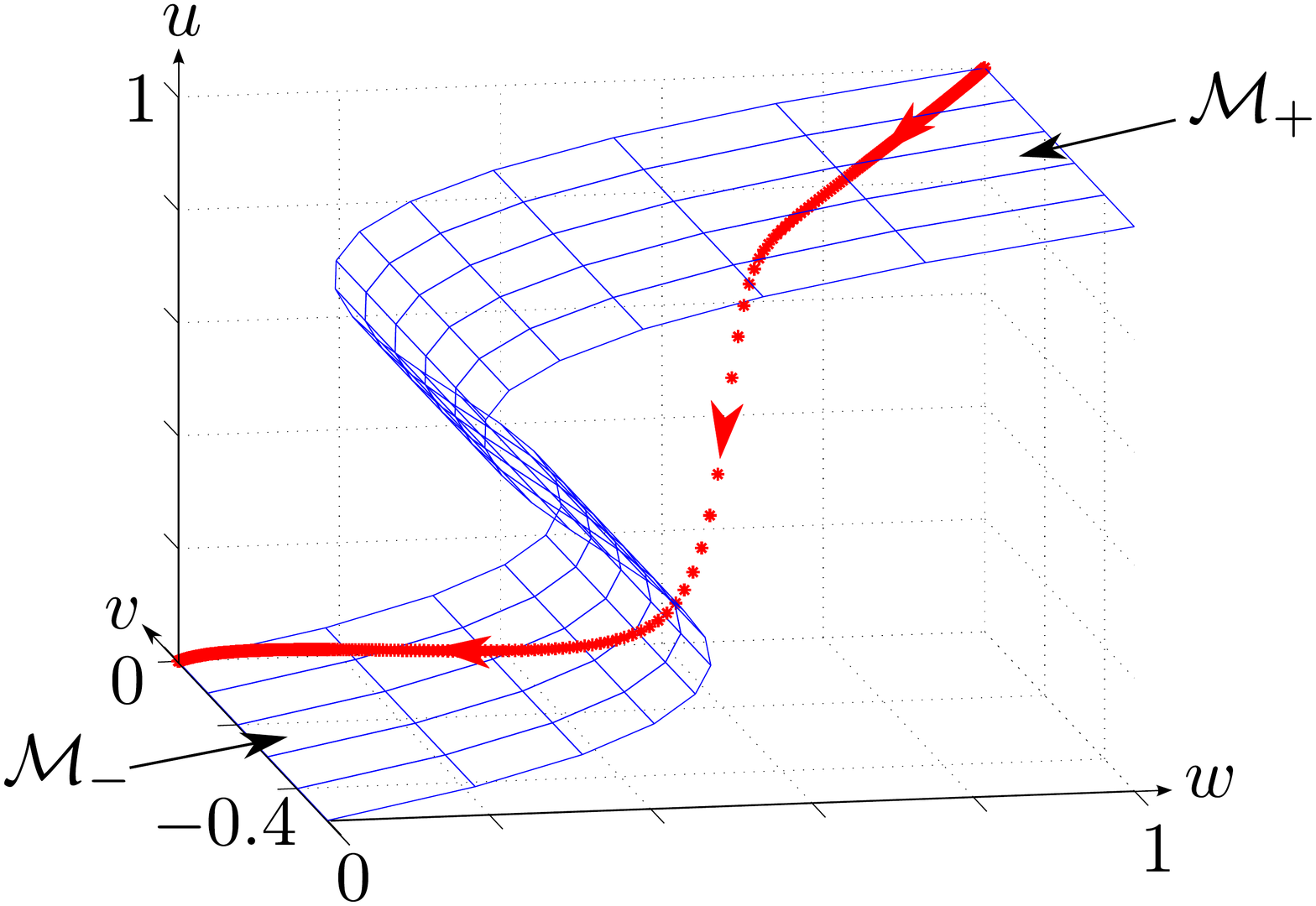}
\caption{Numerically computed heteroclinic ($a=0.4$,  $d=1/8$) in \eqref{eq:sa_3ODE}, together with slow manifolds $\mathcal{M}_\pm$.}
\label{fig:cubicsurface}
\end{figure}

The location of the fold points is of course given just by the extrema of $g$, which we computed earlier as  $u_\pm$; see Figure \ref{fig:vwtrajectory}.

\begin{figure}[h]
  \includegraphics[width=0.5\textwidth]{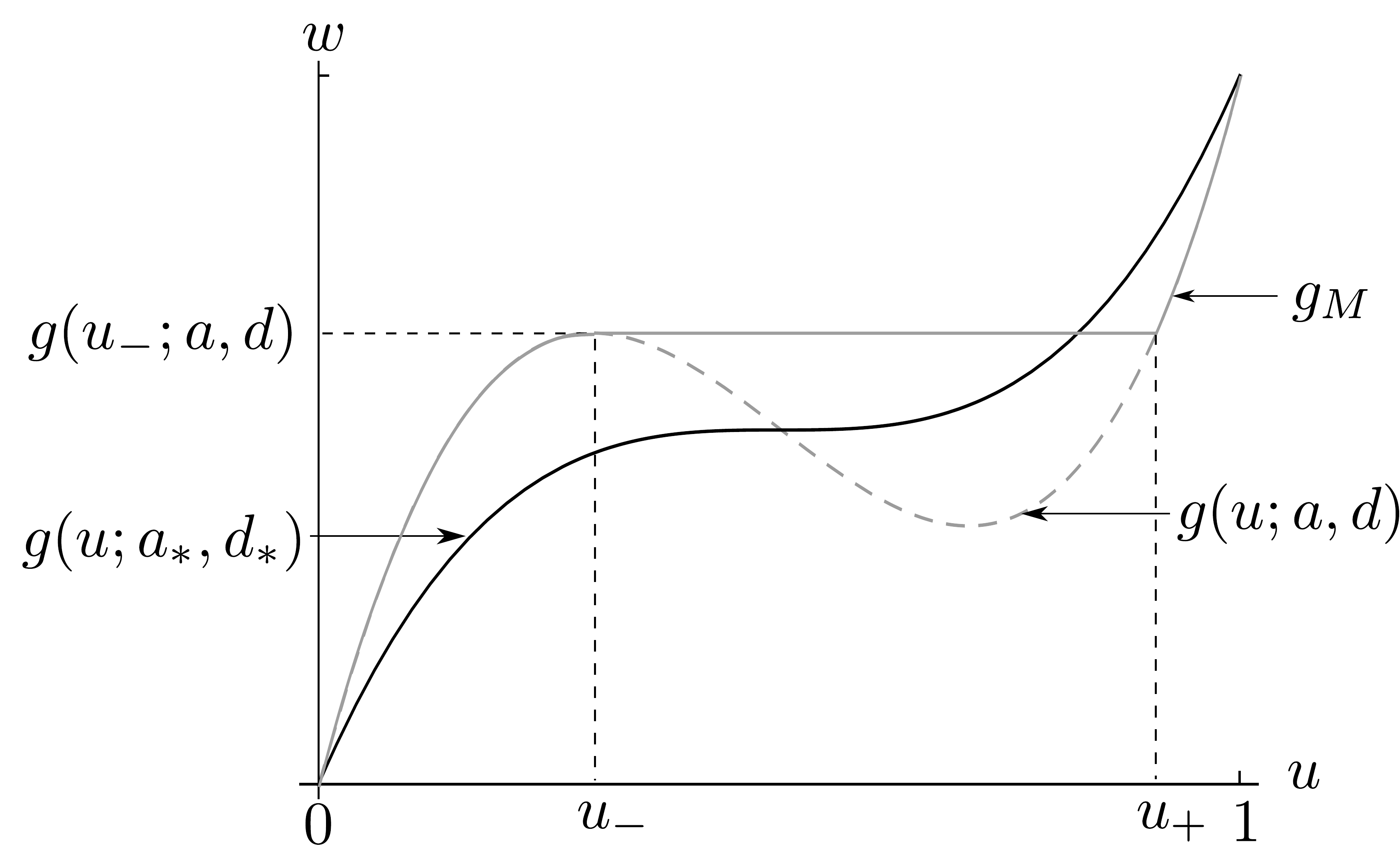}\hfill
  \includegraphics[width=0.49\textwidth]{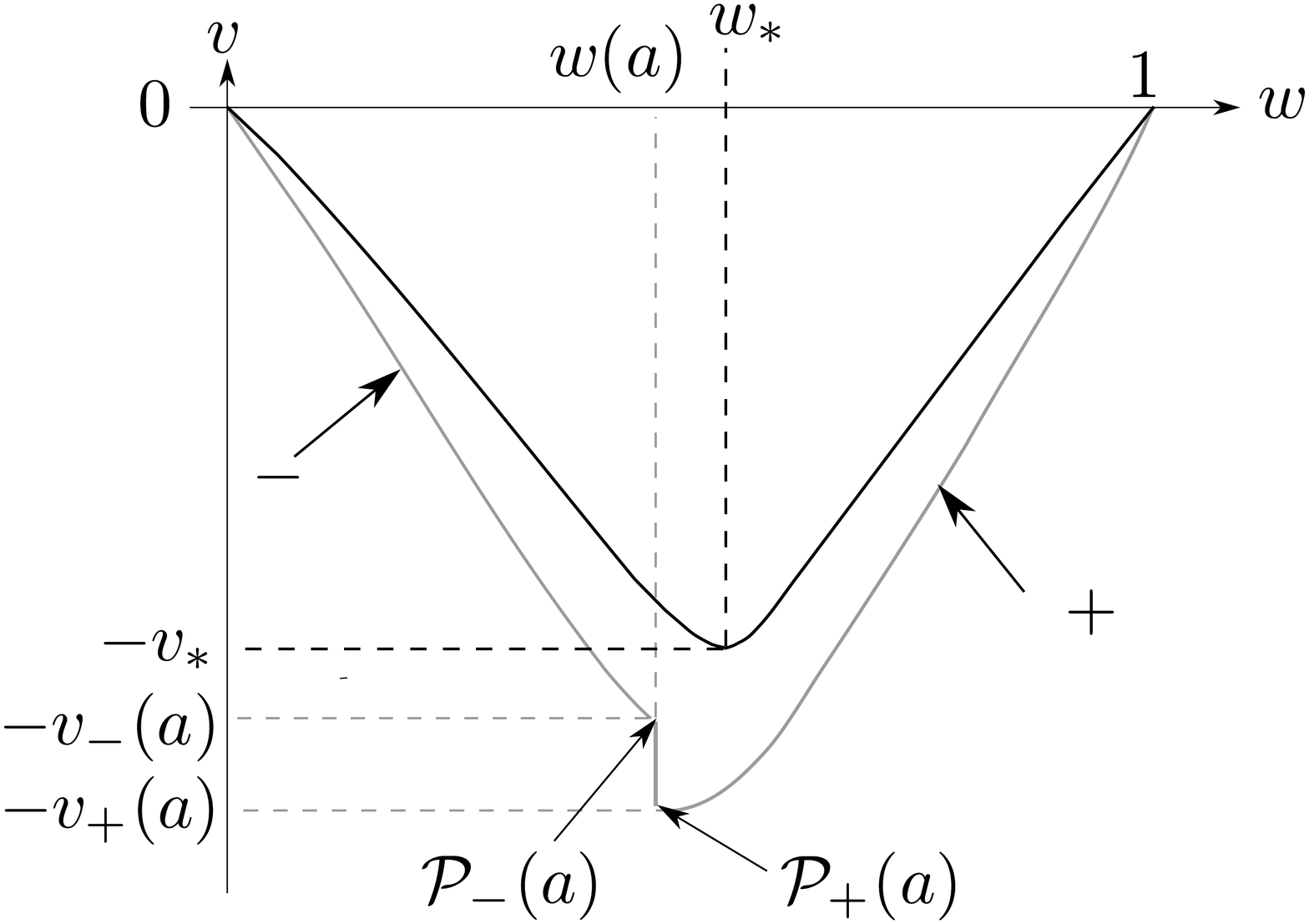}
 \caption{The nonlinearity $g(u;a,d)$ with labeled maxima and fast jump (left). Projection of singular trajectory (black) and nearby trajectory (gray) into the $(w,v)$-plane (right).}
\label{fig:vwtrajectory}
\end{figure}

Since equilibria $\underline{u}_\pm$ lie on different branches, $\mathcal{M}_+$ and $\mathcal{M}_-$, respectively, we are looking for heteroclinics between the two manifolds through the flow of (\ref{e:fast}) at $c=0$. Such heteroclinics, however, exist only at the fold points, where they are simple vertical lines in $\underline{u}$-phase space. They will require a more detailed analysis in the sequel, and are responsible for the peculiar scaling of wave speeds.

\paragraph{Slow system.} 
We next analyze the slow system, setting $c=0$ and formally obtaining an algebro-differential equation. We solve the last equation for $u$ as a function of $w$ and substitute into the remaining two equations, to obtain
\begin{equation}
\begin{split}
w_x & = v,\\
v_x & = w-g_\mp^{-1}(w).
\end{split}\label{e:pend}
\end{equation}
These equations determine the slow flow on the manifolds $\mathcal{M}_\pm$ in the sense of \cite{fenichel:79}. The slow system (\ref{e:pend}) is a Hamiltonian system with conserved Hamiltonian 
\begin{equation}
H_{\mp} (w,v) =\frac{1}{2}v^2+G_{\mp}(w),
\label{eq:g}
\end{equation}
where
\begin{equation}
G_{\mp}(w) =-\frac{1}{2}w^2+\frac{3}{4d}(g^{-1}_{\mp}(w))^4-\frac{2}{3d}(1+a)(g^{-1}_{\mp}(w))^3+\frac{1}{2}\left(1+\frac{a}{d}\right)\left(g^{-1}_{\mp}(w)\right)^2.
\label{eq:h}
\end{equation}
We can therefore determine the location of stable and unstable manifolds of $\underline{u}_-$ and $\underline{u}_+$, respectively, rather explicitly. 

\paragraph{The singular trajectory.} The singular trajectory consists of the (slow) unstable manifold $\mathcal{W}^\mathrm{u}_+$ of $\underline{u}_+$ in $\mathcal{M}_+$, the stable manifold $\mathcal{W}^\mathrm{s}_-$ of  $\underline{u}_-$ in $\mathcal{M}_-$, both connected by a vertical trajectory of the fast flow. In order for such a trajectory to exist (without jump), we need to match stable and unstable manifolds at a fold point. Let $w(a)=g(u_-;a,d)$ be the value of the maximum of $g$ and and let $\mathcal{P}_+(a)$ and $\mathcal{P}_-(a)$ be the points on $\mathcal{W}^\mathrm{u}_+$ and $\mathcal{W}^\mathrm{s}_-$, respectively, intersected with $w=w(a)$, and projected into the $(w,v)$-plane; see Figure \ref{fig:vwtrajectory}. Using the explicit expression for the Hamiltonian (\ref{eq:g}) and (\ref{eq:h}), one finds that $\mathcal{P}_-(a)=\mathcal{P}_+(a)$ precisely at the pinning boundary; that is at $a=a_-$ with $w=w_*$ and $v=-v_*$.

\paragraph{Finite $c>0$ --- persistence of slow manifolds.}
For finite $c$, slow manifolds persist and  $\mathcal{W}^\mathrm{u}_+$ and $\mathcal{W}^\mathrm{s}_-$, respectively, vary smoothly in $c$, up to a neighborhood of the vertical plane $w=w_*$. Also, the stable foliation varies smoothly in $c$ so that we are left with matching the two-dimensional unstable manifold of $\underline{u}_+$ with the one-dimensional stable manifold of $\underline{u}_-$ in a vicinity of the fold point of $\mathcal{M}_-$, near $\mathcal{P}_-(a_-)$. We emphasize that up to a neighborhood of this fold point, manifolds vary smoothly in $c$, so that matching corrections are $\rmO(c)$. The remainder of the proof is based on an analysis of the dynamics near the fold point, and a computation of the $\rmO(c)$-corrections to stable and unstable manifolds outside of this neighborhood. Both are then combined to obtain a reduced matching equation, which gives the desired expansions.

\paragraph{Slow passage through the fold.}
The fast system involves traveling around a fold. This type of system was analyzed in \cite{krupa} and we make use of these results via coordinate transformations. Based on locations of extrema of $g$, we define
\begin{equation*}
w_\ast  := g(u_-) \text{ where } g'(u_-;a_-,d)=0,\qquad
v_\ast  := \sqrt{-2G_-(0)},\qquad
u_\ast  := u_-;
\end{equation*}
see Figure \ref{fig:vwtrajectory}. Now set  $\tilde{w}=w-w_\ast,\: \tilde{v}=v+v_\ast,\: \tilde{u}=u-u_-$, and obtain
\begin{equation}
\begin{split}
& \tilde{w}_y=c (\tilde{v}-v_\ast),\\
& \tilde{v}_y=c (\tilde{w}-\tilde{u})+c(w_\ast-u_-),\\
& \tilde{u}_y=d(-\tilde{w}-w_\ast + g(\tilde{u} +u_-)).\\
\end{split}
\label{eq:fastsystem}
\end{equation}
We expand this system, using the Taylor expansion of $g$ near the maximum where $\partial_u g(u_\ast;a_-,d)=0$. We therefore write   $\gamma:=\partial_ag(u_\ast;a_-,d)\neq 0$ and $\alpha:=\frac{1}{2}\partial_u^2g(u_\ast;a_-,d)\neq 0$, $\mu=a_--a$, and find
\begin{equation}
\begin{split}
& \tilde{w}_y=-c v_\ast+\rmO(c\tilde{v}),\\
& \tilde{v}_y=c(w_\ast-u_-)+\rmO(|c\tilde{v}|+|c\tilde{w}|),\\
& \tilde{u}_y=-d(\tilde{w} +\alpha \tilde{u}^2+\gamma \mu)+\rmO(|\tilde{u}|^3+|\mu||\tilde{u}|).\\
\end{split}
\label{eq:fastsystem2}
\end{equation}
An order-one coordinate transformation, shifting $\tilde{w}$ and rescaling variables  gives at leading order 
\begin{equation*}
\begin{split}
\dot{w}_1&=c_1 ,\\
\dot{v}_1&=c_1,\\
\dot{u}_1&=w_1-u_1^2.
\end{split}\label{e:sf}
\end{equation*}
This  system is the time-reversed version of a system studied by Krupa and Szmolyan \cite{krupa}, extended by the trivial equation for $v_1$. The dynamics of the $(w,u)$-system is depicted in Figure \ref{fig:criticalmanifold}. As the solution $u$ leaves $\mathcal{M}_+$ and falls past the fold point of the function $g$ the distance from the singular solution depends on $c$ ($\Delta^\text{in}$) in a characteristic $2/3$-scaling \cite{krupa},
\begin{equation*}
\Delta w_1=-\Omega_0 c_1^\frac{2}{3}+\mathcal{O}\left(c_1 | \ln(c_1) | \right),
\end{equation*}
where $\Omega_0\approx 2.3381$ is the first positive root of $J_{-\frac{1}{3}}(\frac{2}{3}z^\frac{3}{2})+J_{\frac{1}{3}}(\frac{2}{3}z^\frac{3}{2})$. Inspecting the construction in \cite{krupa}, it is not difficult to see that one can also obtain these scalings in the three-dimensional situation, with 
\begin{equation*}
\Delta v_1=-\Omega_0 c_1^\frac{2}{3}+\mathcal{O}\left(c_1 | \ln(c_1) | \right);
\end{equation*}
see also \cite{gils,SW}. Reversing the scalings, we find 
\begin{equation}
\begin{split}
\Delta\tilde{w}&=-\alpha \Delta w_1=\Omega_0 \left(\frac{v_\ast^2}{d^2{\alpha} }\right)^\frac{1}{3}c^\frac{2}{3}+\mathcal{O}\left(c\, | \ln(c) | \right),\\
\Delta\tilde{v}&=\frac{w_\ast-u_\ast}{v_\ast}\Delta\tilde{w}+\mathcal{O}\left(c\, | \ln(c) | \right).\\
\end{split}
\label{eq:Delta}
\end{equation}
\begin{figure}
\centering
\includegraphics[height=2.4in]{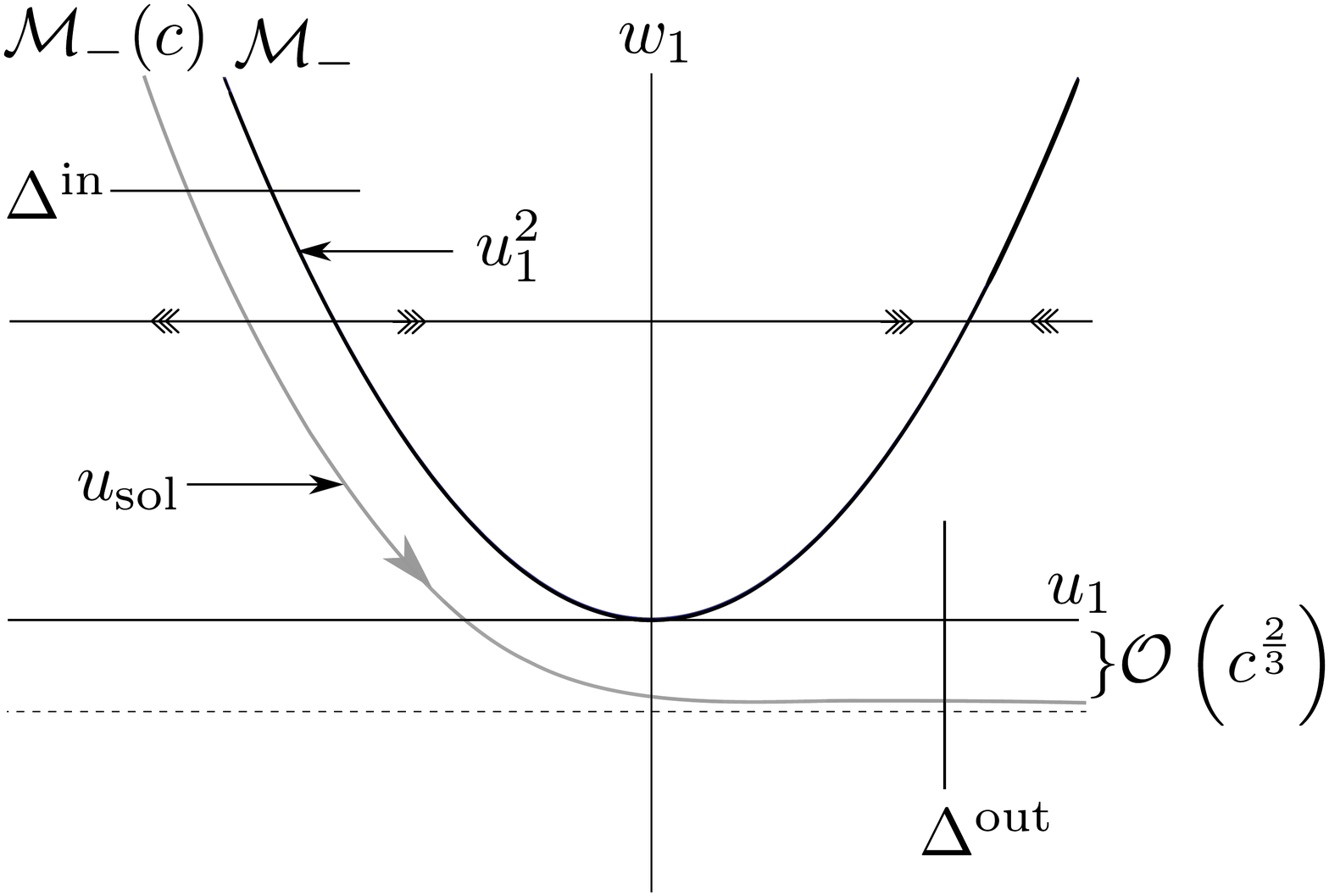}\hfill 
\includegraphics[height=2.4in]{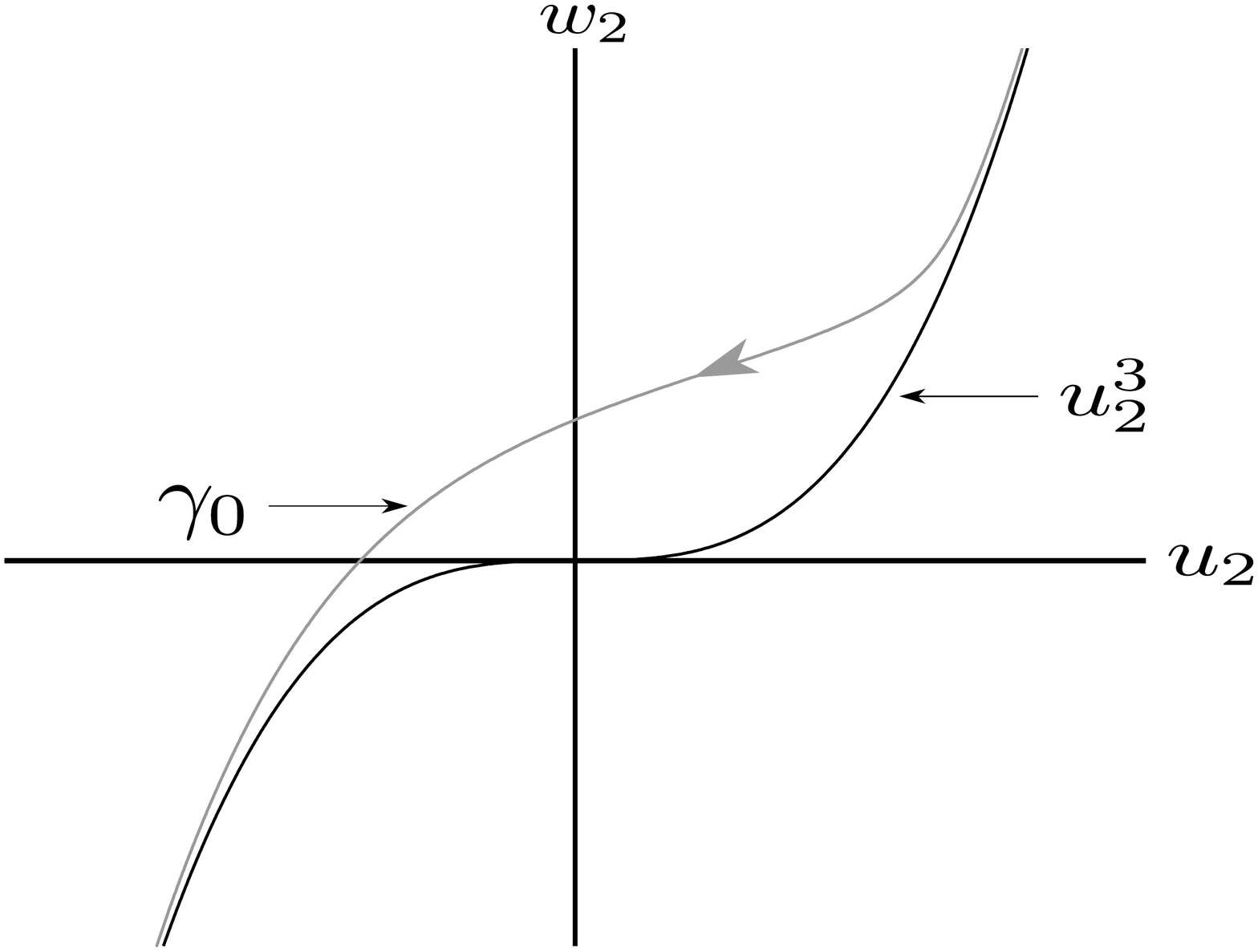}
\caption{Illustration of slow passage through a fold (\ref{e:sf}), left;  \cite{krupa}. Slow passage through an inflection point (right), as used in (\ref{eq:k}).}
\label{fig:criticalmanifold}
\end{figure}

\paragraph{Matching stable and unstable manifolds above the fold.} After the passage through the fold, we have expansions for the location of the stable manifold $\mathcal{W}^\mathrm{s}_+$ in a cross section to the fast flow slightly above the fold, $u=u_*+\delta$, $(w,v)\in\R^2$. In this cross section, $\mathcal{W}^\mathrm{u}_-$ forms a curve that contains $\mathcal{W}^\mathrm{s}_-$ at $c=0$, since we constructed a singular continuous trajectory. We will now need to compute expansions of $\mathcal{W}^\mathrm{u}_+$ in parameters $a$ and $c$ in order to find intersections. Since $\mathcal{W}^\mathrm{u}_+$ is smooth in both variables, the leading order term in $c$ stems solely from the passage through the fold, $c^{2/3}$. It is therefore sufficient to calculate expansions in $a$. Notice that the direction of the fast flow is vertical, independent of $a$, so that $a$-derivatives of the unstable manifold can be calculated from the slow flow on $\mathcal{M}^+$, only, at $c=0$, in (\ref{e:pend}). 

Figure \ref{fig:vwtrajectory} illustrates this matching process. Denote the fold location by $w(a)$ and write $\mathcal{P}_-(a)$ for the intersection of the stable manifold with the fold, before the corrections due to the passage through the fold. On the other hand, $\mathcal{P}_+(a)$ shall denote the intersection of the unstable manifold $\mathcal{W}_+^\mathrm{u}$ with $w=w(a)$. We write $-v_\pm(a)$ for the (negative) $v$-coordinates of $\mathcal{P}_\pm(a)$. 
At $a_-$, we have 
\begin{equation*}
\mathcal{P}_\ast=\mathcal{P}_\pm(a_-)=(w_\ast,-v_\ast).
\end{equation*}
From the Hamiltonian structure, we find $H(\mathcal{P}_-(a))=H((0,0))=0$, so that
\begin{equation}
w(a)  = g(u_-;a,d),\qquad 
v_-(a)=\sqrt{-2G_-(w(a))}.
\label{e:v-}
\end{equation}
Similarly, $H(\mathcal{P}_+(a)) =H((1,0))=\frac{1}{12d}(1+2a)$ so that 
\begin{equation}
w(a)  = g(u_+;a,d),\qquad
v_+(a)  =\sqrt{2\left(\frac{1+2a}{12d}-G_+(w(a))\right)}.
\label{e:v+}
\end{equation}
In addition, we need the tangent space to the two-dimensional unstable manifold $W^\mathrm{u}_-+$ at $\mathcal{P}_+(a_-)$, intersected with the horizontal cross section, which is simply given by the time derivative of the slow flow in $\mathcal{M}_+$, projected on the $(w,v)$-plane, $t_+=(-v_\ast,w_\ast-u_\ast)$.

We then project the equation for the intersection of $\mathcal{W}_+^\mathrm{u}$ and $\mathcal{W}_-^\mathrm{s}$ in the cross section $u\equiv const$ onto $t_+^\perp$ using Lyapunov-Schmidt reduction, to arrive at the reduced equation 
\begin{equation}\label{e:match}
-\langle\partial_a \mathcal{P}_+(a),d^\perp\rangle\mu=-\langle \partial_a \mathcal{P}_-(a),d^\perp\rangle\mu + \langle s,d^\perp\rangle c^\frac{2}{3}+\rmO(\mu^2+| c\log c|),
\end{equation}
where 
\begin{equation}
s:=(\Delta\tilde{w},\frac{w_\ast-u_\ast}{v_\ast}\Delta\tilde{w})=\Omega_0 \alpha\left(\frac{v_*}{\alpha^2 d}\right)^\frac{2}{3} c^\frac{2}{3}(1,\frac{w_\ast-u_\ast}{v_\ast}).
\end{equation}
Solving (\ref{e:match}) for $\mu$, we find that:
\[\mu=\frac{\langle s,d^\perp\rangle}{\langle\partial_a\mathcal{P}_-(a),d^\perp\rangle - \langle\partial_a\mathcal{P}_+(a),d^\perp\rangle}
 c^\frac{2}{3}+\mathcal{O}\left( c\, | \ln( c) | \right).\]
and, after substitution, 
\begin{equation*}
\mu=\frac{\Omega_0 \alpha \left(\frac{v_\ast}{\alpha^2 d}\right)^\frac{2}{3}(w_\ast-u_\ast)}{\partial_a v_+(a) 
- \partial_a v_-(a)} c^\frac{2}{3}+\mathcal{O}\left( c\, | \ln( c) | \right),
\end{equation*}
where $v_\pm$ where defined in \eqref{e:v-},\eqref{e:v+}, and $G_\pm$ were defined in \eqref{eq:h}. In particular, we find, using the notation from Theorem \ref{thm:2D},
\begin{equation}
k_1=\left(\frac{\partial_a v_+(a_-) 
- \partial_a v_-(a_-)}{\Omega_0 \alpha \left(\frac{v_\ast}{\alpha^2 d}\right)^\frac{2}{3}(w_\ast-u_\ast)}\right)^{\frac{3}{2}}.
\label{eq:k1}
\end{equation}

\end{Proof}

\subsection{Asymptotics at the tip of the pinning region}
\label{subsec:prtip}

Inspecting the shape of the pinning region, Figure \ref{fig:generalcubicpinningregion}, we expect a transition at $d=1/4$. The previous result, Theorem \ref{thm:2D}, describes speeds for $a$ close to the left boundary, when $d<1/4$. For $d>1/4$, interfaces are stationary at $a=1/2$, only, and $c\sim \mu$ is smooth. It is therefore interesting to examine speeds at criticality, fixing $d=1/4$ and varying $\mu=a-1/2$ near the origin. We will see below that a very similar strategy leads to expansions with a new exponent, $5/4$. The heart of the analysis, however, relies on a singular perturbation problem that involves the slow passage through an inflection point, which has not been studied in a rigorous fashion, to our knowledge. We outline the key elements of such a study, but do not develop a full geometric proof. For a rigorous geometric approach of the related slow passage through a cusp, we refer to the recent study in \cite{broer}; the results there cover a more general unfolding but do not give expansions 
for our case. 

Our result from this not fully rigorous study is formulated as follows.

\begin{main}
For fixed $d=d_\ast$, as $a$ approaches the boundary of the pinning region, $a\nearrow a_\ast=1/2$, the wave speed of the interface from Theorem \ref{thm:bates} is
\begin{equation*}
c=k_c(a_*-a)^\frac{5}{4}+\rmo((a_*-a)^{\frac{5}{4}}).
\end{equation*}
The constant $k_c$ is given ``explicitly'' in equation \eqref{eq:kc}. 
\label{thm:fivefourths}
\end{main}

In order to arrive at the asymptotic expansion, we follow the strategy of the proof of Theorem \ref{thm:2D}. 

\paragraph{Fast system}
Again, we rescale time and in shifted coordinates $(\tilde{w}, \tilde{v}, \tilde{u})$ defined as in \eqref{eq:fastsystem}, where now $w_*,-v_*,u_-$ denote the inflection point of the slow manifold,
\begin{equation*}
\begin{split}
\tilde{w}_y&=c v_*+c \tilde{v},\\
\tilde{v}_y&=c(w_*-u_*)+c(\tilde{w}-\tilde{u}),\\
\tilde{u}_y&=-d_*(\tilde{w}+w_*-g(\tilde{u}+u_*)).
\label{eq:i}
\end{split}
\end{equation*}
Expanding $g$, we find at leading order
\[
\tilde{u}_y=-d_*(\tilde{w}+\mu-4\tilde{u}^3-4\tilde{u}^2\mu).
\]
We shift $\tilde{w}$ to eliminate the linear term in $\mu$, to obtain
\[
\tilde{u}_y=\tilde{u}^3-\frac{1}{4}\tilde{w}.
\]
Since we are situated at the critical values $w_*=u_*=a_*=\frac{1}{2}$ we see that the equation for $\tilde{v}$ simplifies to
\[
\tilde{v}_y =c(\tilde{w}-\tilde{u}).
\]
Altogether, substituting $d_*=1/4$, we find at leading order 
\begin{equation}
\begin{split}
\tilde{w}_y&=c v_*,\\
\tilde{v}_y&=c(\tilde{w}-\tilde{u}),\\
\tilde{u}_y&=\tilde{u}^3-\frac{1}{4}\tilde{w}.
\end{split}
\label{eq:j}
\end{equation}
In analogy to the slow passage through the fold discussed in \cite{krupa}, we focus on the crucial dynamics in the rescaling chart. While the quadratic in the fold yields a quadratic nonlinearity in a Riccati equation, we obtain here, not surprisingly, a cubic nonlinearity. To be specific, rescale first $\bar{w}=\frac{\tilde{w}}{4}$ and $c_2=\frac{c v_*}{4}$, then set $\bar{w}=c_2^\frac{3}{5}w_2, \tilde{u}=c_2^\frac{1}{5}u_2,
\tilde{v}=c_2^\frac{4}{5}v_2,$ and (rescaling time) $y=c_2^\frac{-2}{5}z,$ to obtain (denoting $z$-derivatives with $\dot{\ }$),
\begin{equation*}
\begin{split}
\dot{w_2}&=-1+c_2^\frac{4}{5}\frac{v_2}{v_*},\\
\dot{v_2}&=-\frac{4}{v_*}u_2+\frac{4^2}{v_*}c_2^\frac{2}{5}w_2,\\
\dot{u_2}&=u_2^3-w_2.
\end{split}
\end{equation*}
At $c=0$, we finally obtain
\begin{equation}
\begin{split}
\dot{w_2}&=-1,\\
\dot{v_2}&=\frac{4}{v_*}u_2,\\
\dot{u_2}&=u_2^3-w_2.\\
\label{eq:k}
\end{split}
\end{equation}
We will show in the appendix that the following facts hold for \eqref{eq:k}:
\begin{enumerate}
\item There exists a unique trajectory  $\gamma_0$ in the $(w_2,u_2)$-plane such that $u_2(z)^3-w_2(z)\to 0$ for $z\to \pm \infty$. 
%This trajectory can be written in the form  $\gamma_0=\{(h_0(w),w)\,|\,w \in \R\}$, with $h_0$ monotone, and  $h_0(w) \sim w^\frac{1}{3}$ as $w \rightarrow \pm \infty$. 
\item The Cauchy Principal Value $C_0$ of $u_2$ exists and
$$C_0 := P.V. \int_\R u_2(z)\rmd z<0.$$
%\item There exists a unique trajectory  $\gamma_0$ in the $(u_2, w_2)$ plane such
%that $u_2(z)^3-w_2(z)\to 0$ for $z\to \infty$. This trajectory can be written in the form  $\gamma_0=\{(h_0(w),w)\,|\,w \in \R\}$, with $h_0$ monotone, and  $h_0(w) \sim w^\frac{1}{3}$ as $w \rightarrow \pm \infty$; see Figure \ref{fig:criticalmanifold} (right panel). 
%\item The Cauchy Principal Value $C_0$ of $h_0(w_2)$ exists and
%$$C_0 := P.V. \int_\R h_0(w)\rmd w<0.$$
\end{enumerate}
We numerically evaluated the constant as
\[
C_0\approx -2.6524.
\]
From this, we can calculate the ``jump'' of the $v$-component after passage through the inflection point, 
\begin{equation*}
\lim_{T \to \infty} \left(v(T)-v(-T)\right)=C_0,
\end{equation*}
and, in the original variables,  $s_2=\Delta \tilde{v}=C_0(\frac{4}{v_*})^\frac{1}{5}c^\frac{4}{5}+ \rmo(c^\frac{4}{5}).$ 

\paragraph{Matching stable and unstable manifolds.}
Having calculate the effective correction in $v$ after passage through the inflection point, we can match with $a$-derivatives, in a procedure completely equivalent to the fold case. We find $t_+=(1,0)$, so that only the second component $s_2$ of the jump $s$ in the $v$-component contributes. After a calculation analogous to the fold case, we find
\[\mu=\frac{-C_0(\frac{4}{v_*})^\frac{1}{5}}{\langle\partial_a\mathcal{P}_-(a),d^\perp\rangle - \langle\partial_a\mathcal{P}_+(a),d^\perp\rangle}c^\frac{4}{5}+\rmo(c^\frac{4}{5}).\]
Again,  $v_\pm$ where defined in \eqref{e:v-},\eqref{e:v+}, and $G_\pm$ were defined in \eqref{eq:h}. 
In particular, we have
\begin{equation}
k_c=\left(\frac{\partial_av_+(a)-\partial_av_-(a)  }{C_0(\frac{4}{v_*})^\frac{1}{5}}\right)^{\frac{5}{4}}.
\label{eq:kc}
\end{equation}

\subsection{One-sided, first-order kernel}
\label{subsec:onederivative}
A second, interesting example, where asymptotics can be described explicitly, is the one-sided kernel with Fourier symbol $(1+\rmi\ell)^{-1}$, already discussed in Section \ref{subsec:1dpin}. 
Moving interfaces solve 
\begin{equation}
\begin{split}
w_x&=u-w,\\
-c u_x&=d(-u+w)+f_a(u).\\
\end{split}
\label{eq:sa_2ODE1}
\end{equation}
Again, we are interested in wave speed asymptotics near the boundary of the pinning region. The analysis and the result are quite similar to, in fact simpler than in the case of the second order kernel considered in Section \ref{subsec:speedsecondderivative}.
\begin{theorem}
Consider interfaces connecting $u=1$ at $x=-\infty$ to $u=0$ at $x=+\infty$. 
For fixed $d>0$, as $a$ approaches the boundary of the pinning region at $a_\ast$, there exists a unique traveling wave with speed $c$, where 
\begin{equation*}
c=k_2(a_\ast-a)^{\frac{3}{2}}+\mathcal{O}\left((a_\ast-a)^2\ln(a_\ast-a)\right), \text{ as } a\nearrow a_*.
\end{equation*}
The constant $k_2$ is given explicitly in \eqref{eq:k2}.
\label{thm:1D}
\end{theorem}
\begin{Proof}
We proceed in analogy to Section \ref{subsec:speedsecondderivative}. 
Figure \ref{fig:singulartrajectory} contains a phase portrait with slow manifold $w=g(u;a,d)$ and trajectories passing near the fold point of the slow manifold. We are seeking heteroclinic orbits in the $(w,u)$-plane connecting $(1,1)$ and $(0,0)$. For $c<0$, both equilibria are repellers and a heteroclinic cannot exist. For $c>0$, $(1,1)$ is a saddle with one-dimensional stable manifold given through the fast fibration. The origin possesses a one-dimensional unstable manifold which coincides with the slow manifold. The singular trajectory exists when the value of $g$ at the local maximum $u_-$  is 1, precisely at the boundary of the pinning region when $d=\frac{a^2}{4}$. 
\begin{figure}
\centering
\includegraphics[width=0.5\textwidth]{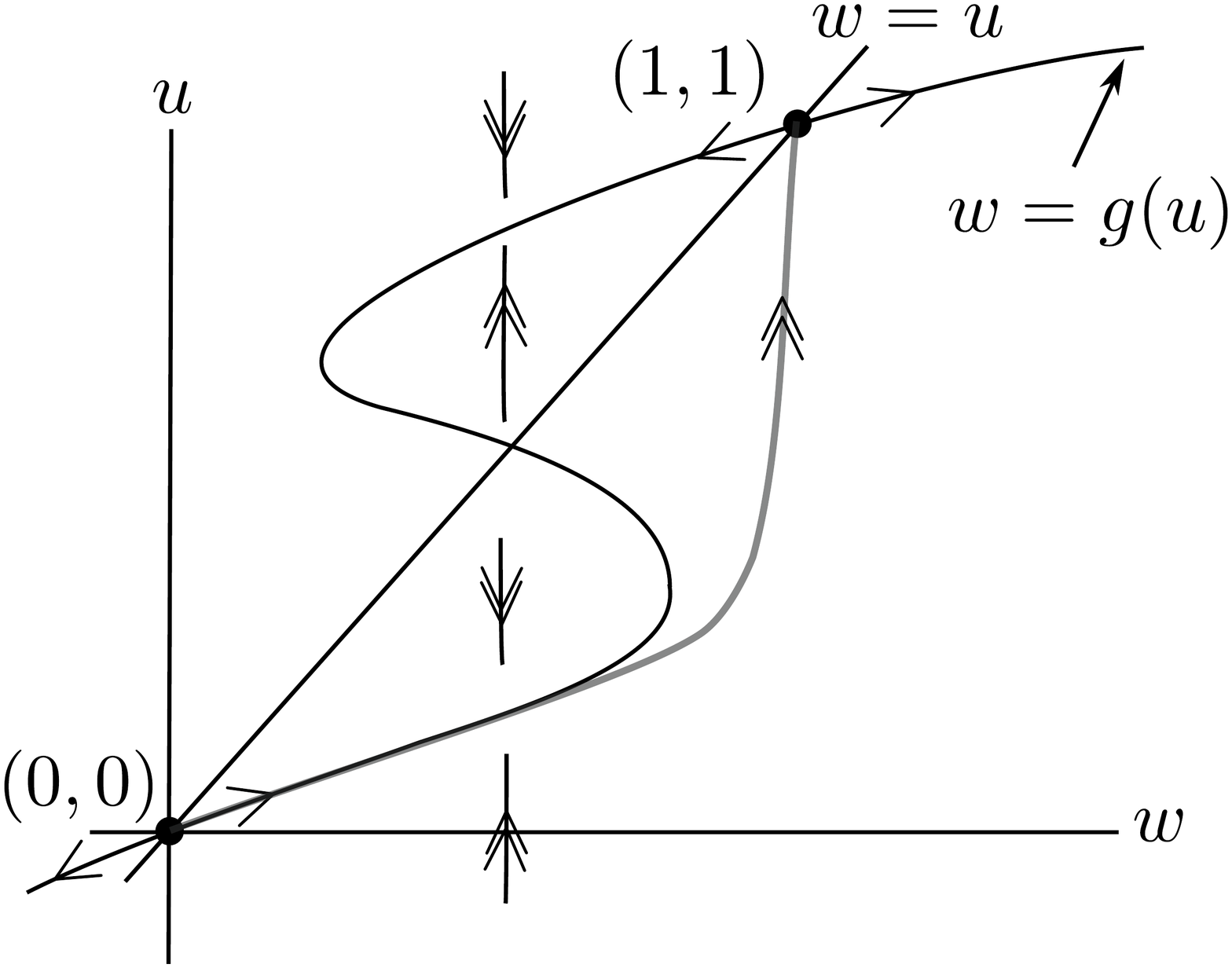}
\caption{Heteroclinic connection and slow-fast structure of \eqref{eq:sa_2ODE1}.}
\label{fig:singulartrajectory}
\end{figure}
When $a$ is slightly decreased, $g(u_-,a,d_\ast)$ is slightly less than $1$. The resulting mismatch between fast fibration and fold point is compensated by the correction from the slow passage through the fold. In order to obtain explicit expansions, we again center our system around the jump point. We use coordinate transformations nearly identical to those used before, to obtain at leading order 
\begin{equation}
\begin{split}
\dot{u}_1&=w_1-u_1^2,\\
\dot{w}_1&=c_1,\\
\end{split}
\label{eq:reduced1d}
\end{equation}
with scalings $c_1=\frac{w_\ast-u_-}{d\alpha^2}c$, and $w_1=-\frac{\tilde{w}}{\alpha}$. Here, again $\alpha:=\frac{1}{2}\partial_u^2g(u_\ast;a_\ast,d)\neq 0$.

Now we have to match points across the fast manifold jump. We write $w(a)$ for the value of the maximum of $g$, and find the matching equation by exploiting the expansion for the transition through the fold \eqref{eq:Delta},
\begin{equation}
\partial_a w(a_*)\mu=\Omega_0\alpha\left(\frac{w_\ast-u_-}{d\alpha^2}\right)^\frac{2}{3}c^\frac{2}{3}+\mathcal{O}\left(|c\,\ln(c)|+\mu^2\right).
 \end{equation}
This readily gives the desired expansion 
\begin{equation*}
c=k_2\mu^\frac{3}{2}+\rmO(\mu^2\ln(\mu)),
\end{equation*}
where $k_2$ is given explicitly through 
\begin{equation}
 k_2=\frac{\left[\partial_a w(a_\ast)\right]^\frac{3}{2}d\alpha^2}{(w_\ast-u_-)(\Omega_0\alpha)^\frac{3}{2}}.
\label{eq:k2}
\end{equation}

\end{Proof}

\section{Rational kernels and beyond --- numerical explorations}
\label{sec:numerics}

We explored scaling laws for unpinning of interfaces numerically. We first report on numerical explorations of the validity range of our predictions from Section \ref{sec:asy}. We then show some evidence of universality for the $3/2$ asymptotics in the case of smooth kernels. Finally, we investigate a family of kernels with increasing singularity and study the dependence of scaling exponents on regularity of the kernel. 

\paragraph{Rational kernels --- validity range of main results.}

\begin{figure}[htp]
\centering
\subfigure[$d=0.1$]{
\label{fig:autoD01}
\includegraphics[width=0.45\textwidth]{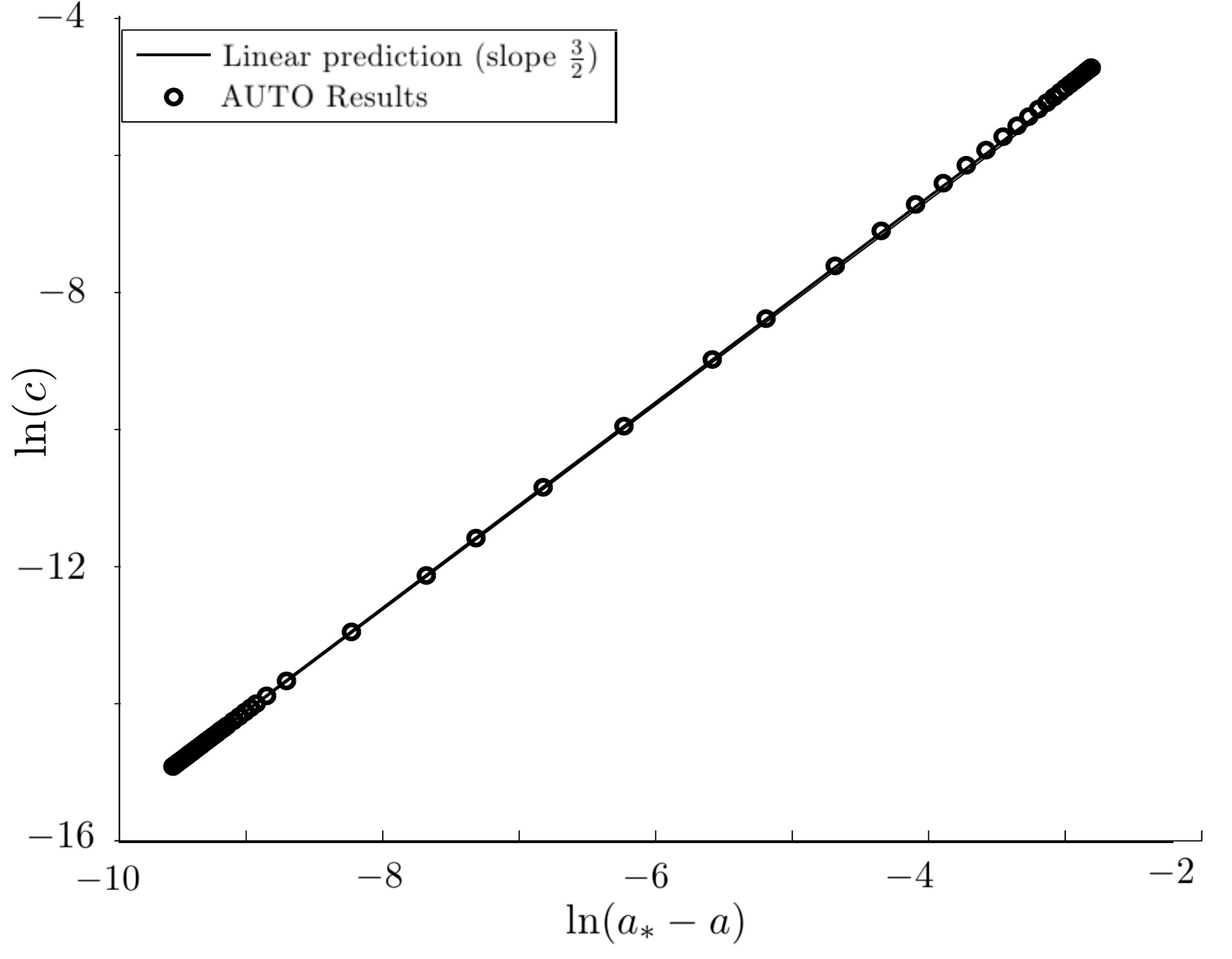}}
\hspace{.3in}
\subfigure[$d=0.1$]{
\label{fig:normKerneld110}
\includegraphics[width=0.45\textwidth]{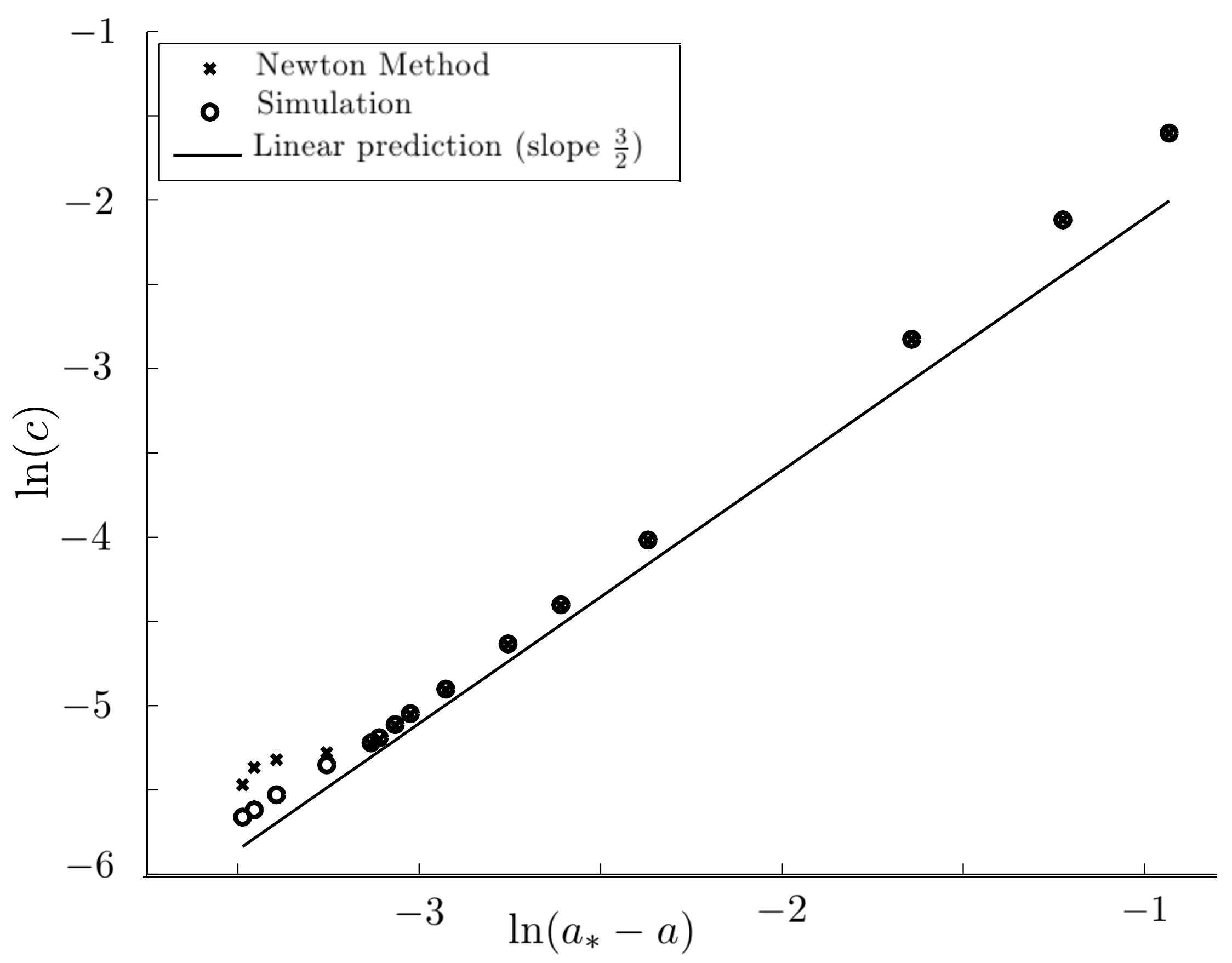}}\\
%bb=0 0 376 326
\subfigure[$d=\frac{1}{4}$]{
\label{fig:autoD25}
\includegraphics[width=0.45\textwidth]{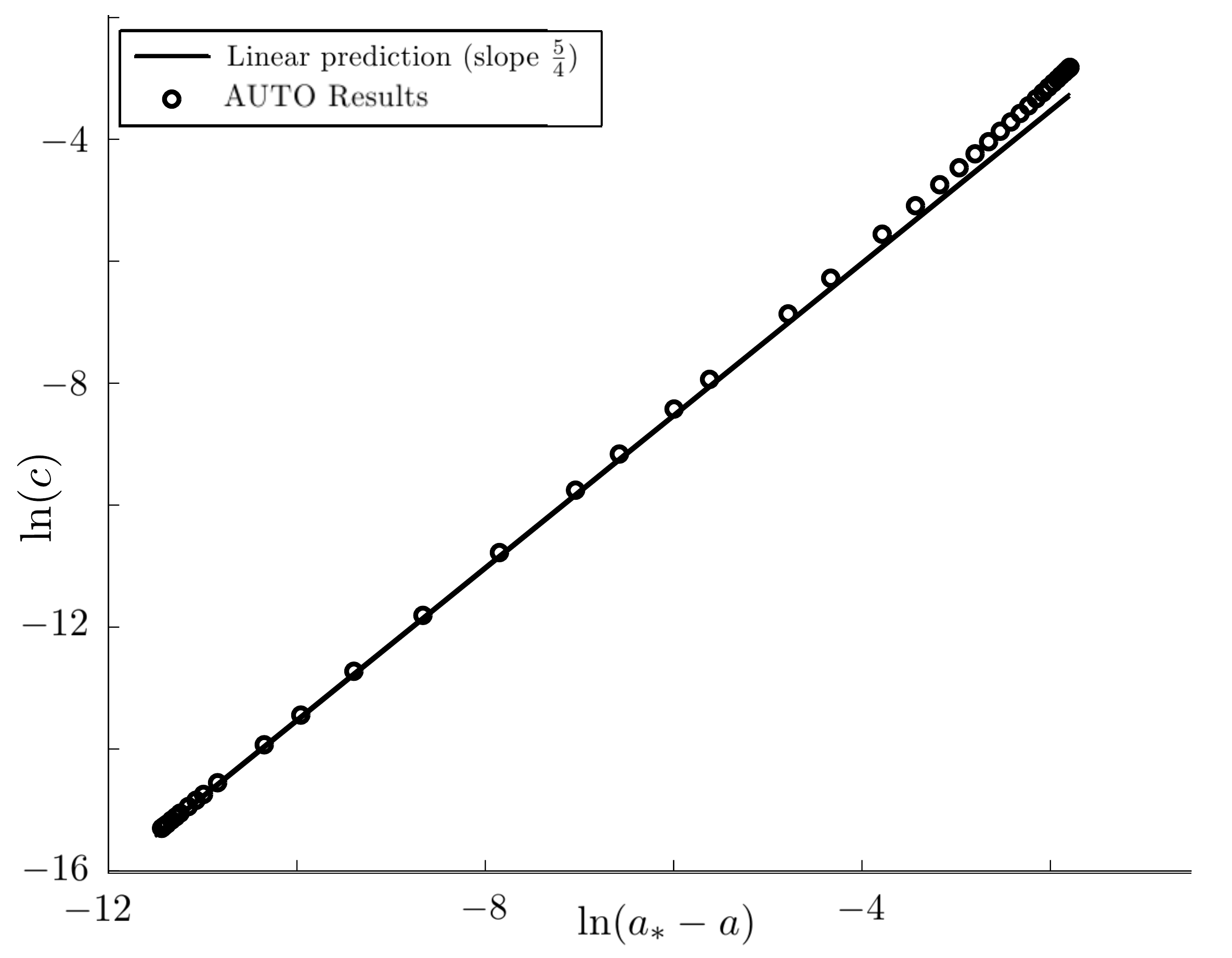}}
%,bb=0 0 398 326
\hspace{.3in}
\subfigure[$d=1$]{
\label{fig:normKerneld1}
\includegraphics[width=0.45\textwidth]{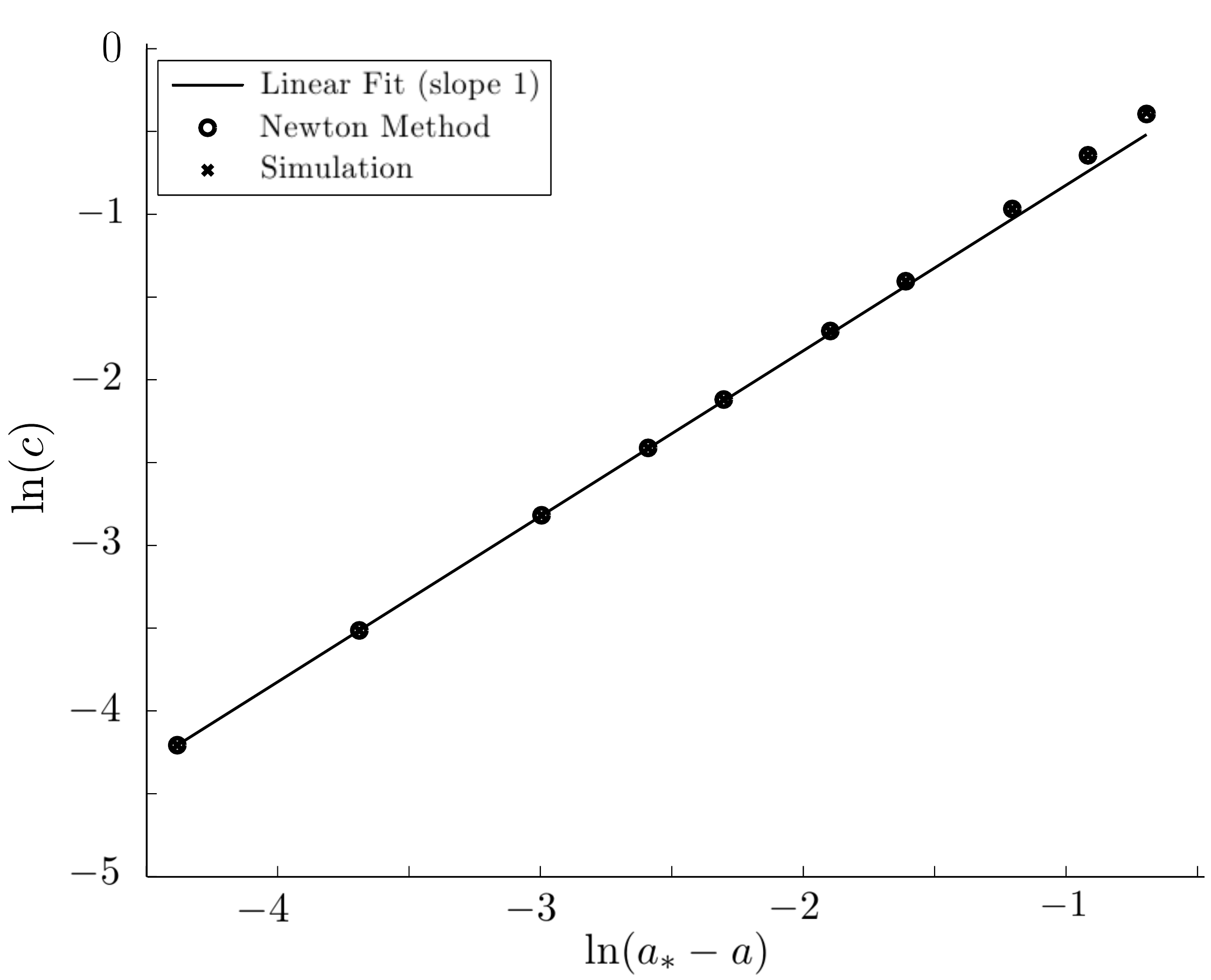}}
%bb=0 0 405 310
\caption{Log-log plots of speed versus distance to pinning boundary. Plotted are results from theoretical predictions, computations with \textsc{auto07p}, and results from 2nd order finite difference approximations using either direct simulation or Newton's method. Best results are, not surprisingly, obtained using \textsc{auto07p}. All results are for $\mathcal{K}(x)=\frac{\mathrm{e}^{-|x|}}{2}$ with $a_*$ given in  Figure \ref{fig:generalcubicpinningregion}.}
\label{fig:numerics}
\end{figure}

Our main results give expansions for the wave speeds near the pinning region. We computed those speeds numerically and compared the scalings. Since $\mu^2\log\mu$-corrections in the expansions are a priori not much smaller compared to leading order terms $\mu^{3/2}$, for moderately small values of $\mu$, one would not necessarily expect strong evidence of asymptotics in direct simulations. All the more surprising, our numerical results show that scalings give in fact excellent predictions for fairly moderate values of $\mu$ and $c$. We computed wave speeds using three methods:
\begin{enumerate}
\item direct simulations;
\item Newton's method for the discretization;
\item \textsc{auto07p}.
\end{enumerate}
In the first two cases, we used second order finite differences in space and simple Euler time stepping. Convolution kernels are evaluated by solving the discretized ODE formulation for $w$, also used in the analysis. In direct simulations, speeds are measured using direct measurements of $\frac{\Delta x}{\Delta t}$, interpolating linearly between grid points. The results are shown in Figure \ref{fig:numerics}. Best results are, not surprisingly, obtained using \textsc{auto07p}. In this case, exponents \emph{and} constants agree with the prediction within roughly 2\%-margins. In detail, our measured\footnote{Discretization details: $dt=0.17$, $dx=0.007324$, on $|x|\leq 30$; best-fit slopes using 4th to 8th smallest data points.}
 and predicted slopes $s_{m/p}$ and $k_{m/p}$ are presented in Table \ref{table:data1}. 

%(a) $k_m=0.5521$, $k_p=0.5457$, $s_m=1.502$, $s_p=1.5$ (b) $k_m=.6842$, $k_p=0.5427$, $s_m=1.544$, $s_p=1.5$ (c) $k_m=0.4030$, $k_p=0.3853$, $s_m=1.259$, $s_p=1.25$ (d) $s_m=1.041$, $s_p=1$.\footnote{Discretization details: $dt=0.17$, $dx=0.007324$, on $|x|\leq 30$; best-fit slopes using 4th to 8th smallest data points.}

\begin{table}[!h] 
\begin{center}
\begin{tabular}{|c||c|c|c|c|c|c|}
\hline & $\mathcal{K}(x)$ & $d$ & $k_m$ & $k_p$ & $s_m$ & $s_p$   \\\hline
Figure \ref{fig:numerics}(a) & $\rme^{-|x|}/2$ & $1/10$ & $0.5521$ & $0.5457$ & $1.502$ & $1.5$ \\\hline
Figure \ref{fig:numerics}(b) & $\rme^{-|x|}/2$ & $1/10$ & $0.6842$ & $0.5427$ & $1.544$ & $1.5$ \\\hline
Figure \ref{fig:numerics}(c) & $\rme^{-|x|}/2$ & $1/4$ & $0.4030$ & $0.3559$ & $1.259$ & $1.25$ \\ \hline
Figure \ref{fig:numerics}(d) & $\rme^{-|x|}/2$ & $1$ & NA & NA & $1.041$ & $1$\\\hline
Figure \ref{fig:numerics2}(a) & $\rme^{-x}\chi_{[0,\infty)}(x)$  & $1/16$ & $0.3892$ & $0.4218$ & $1.492$ & $1.5$  \\\hline
\end{tabular}
\end{center}
\caption{Measured and predicted slopes $s_{m/p}$ and $k_{m/p}$ from Figures \ref{fig:numerics}(a)-(d) and \ref{fig:numerics2}(a).}
\label{table:data1} 
\end{table}

Even in the critical case, we found surprisingly good agreement between predictions and numerical calculations (constant correct to 10\%). Notably, scaling exponents agree well with in the case of (coarse) time stepping, while prefactors differ as expected. 

\begin{figure}[htp]
\centering
\subfigure[$\mathcal{K}=(1+\partial_x)^{-1}$, $d=\frac{1}{16}$]{
\label{fig:auto1D_d0625}
\includegraphics[width=0.45\textwidth]{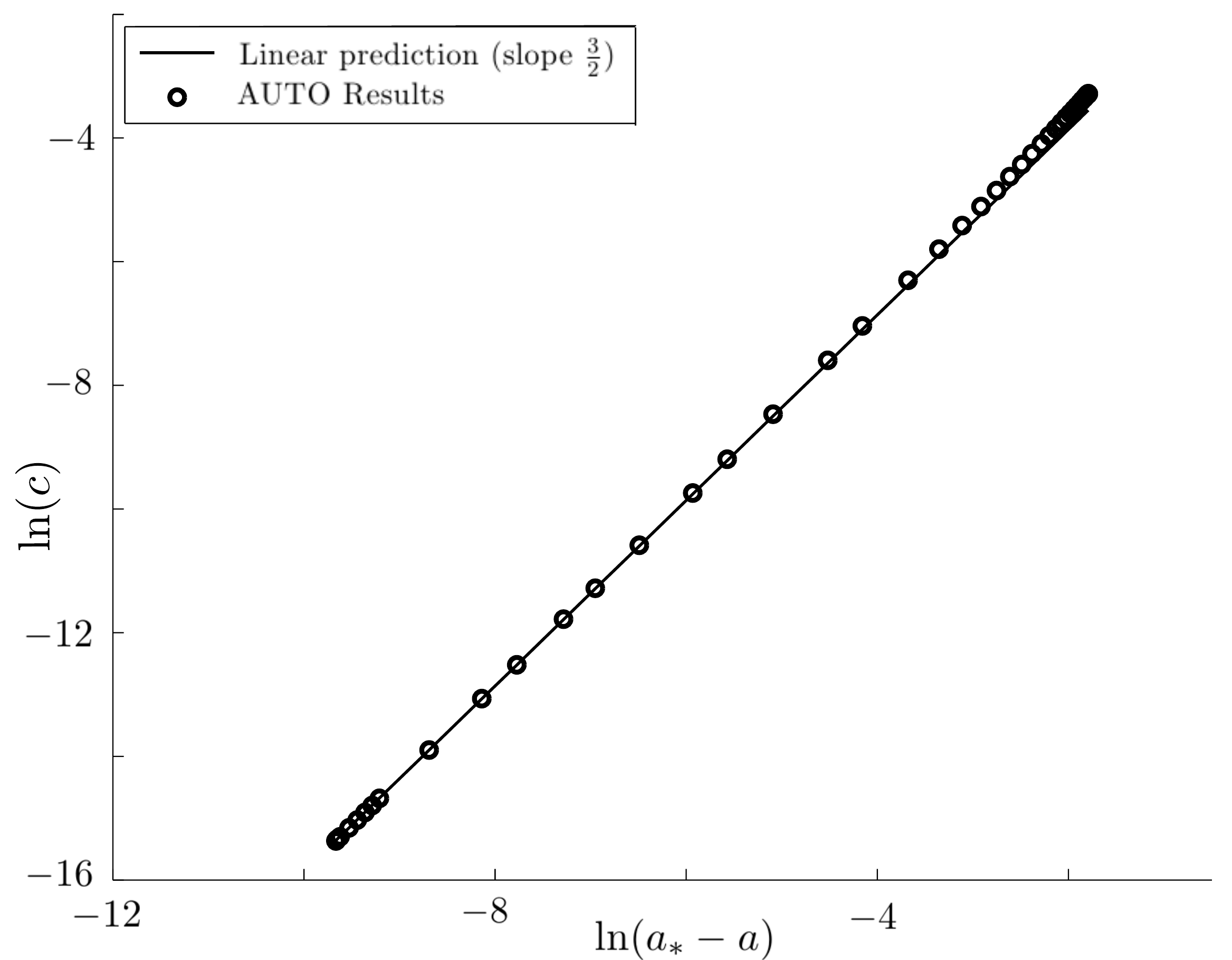}}
\hspace{.3in}
\subfigure[$\mathcal{K}=(1-\partial_x^2 + \partial_4^4)^{-1}$, $d=0.1$]{
\label{fig:2-4DKern}
\includegraphics[width=0.45\textwidth]{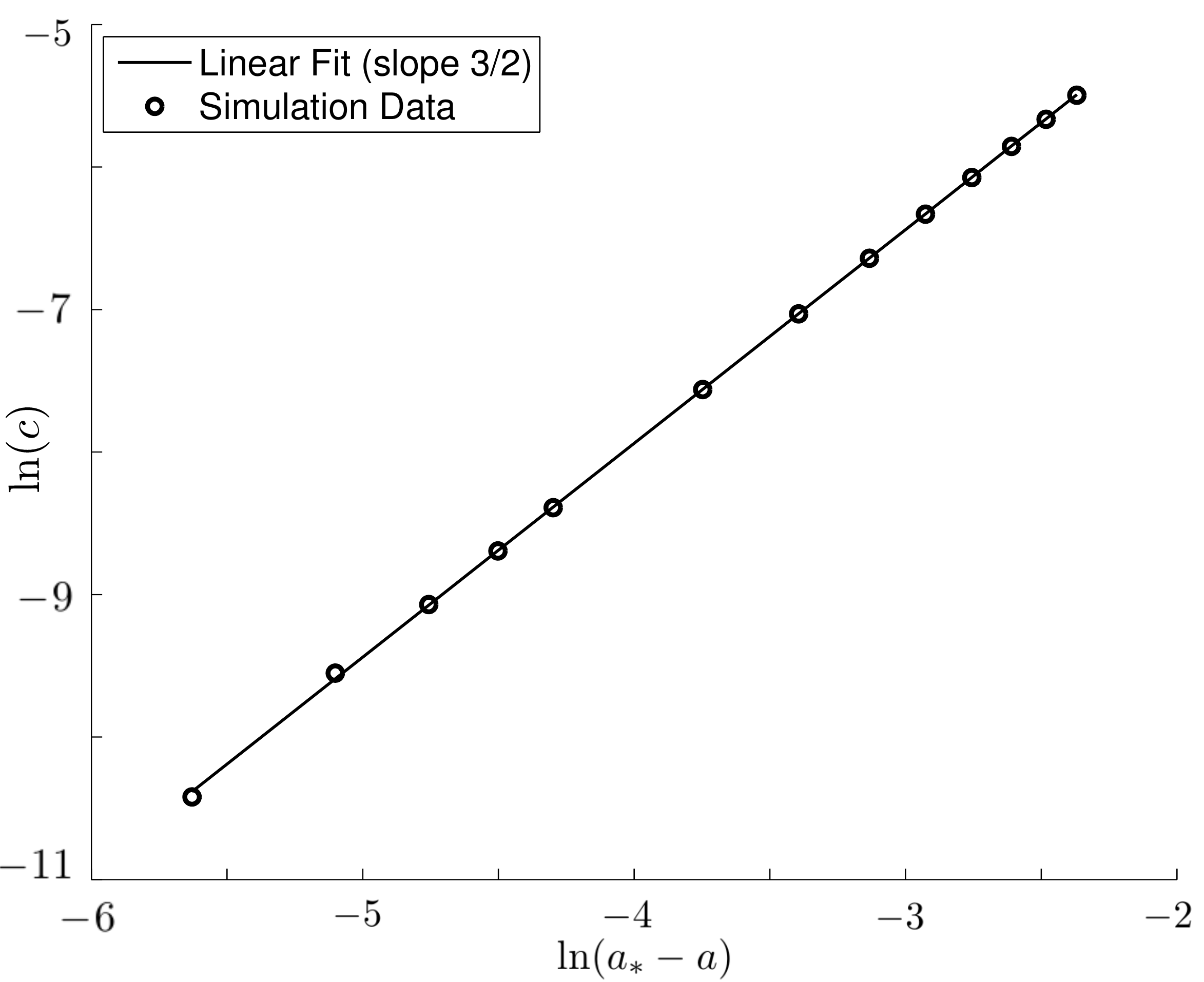}}\\
%bb=0 0 376 326
\subfigure[$\mathcal{K}=\frac{\chi_{[-1,1]}}{2}$, $d=0.1$]{
\label{fig:gaussConv}
\includegraphics[width=0.45\textwidth]{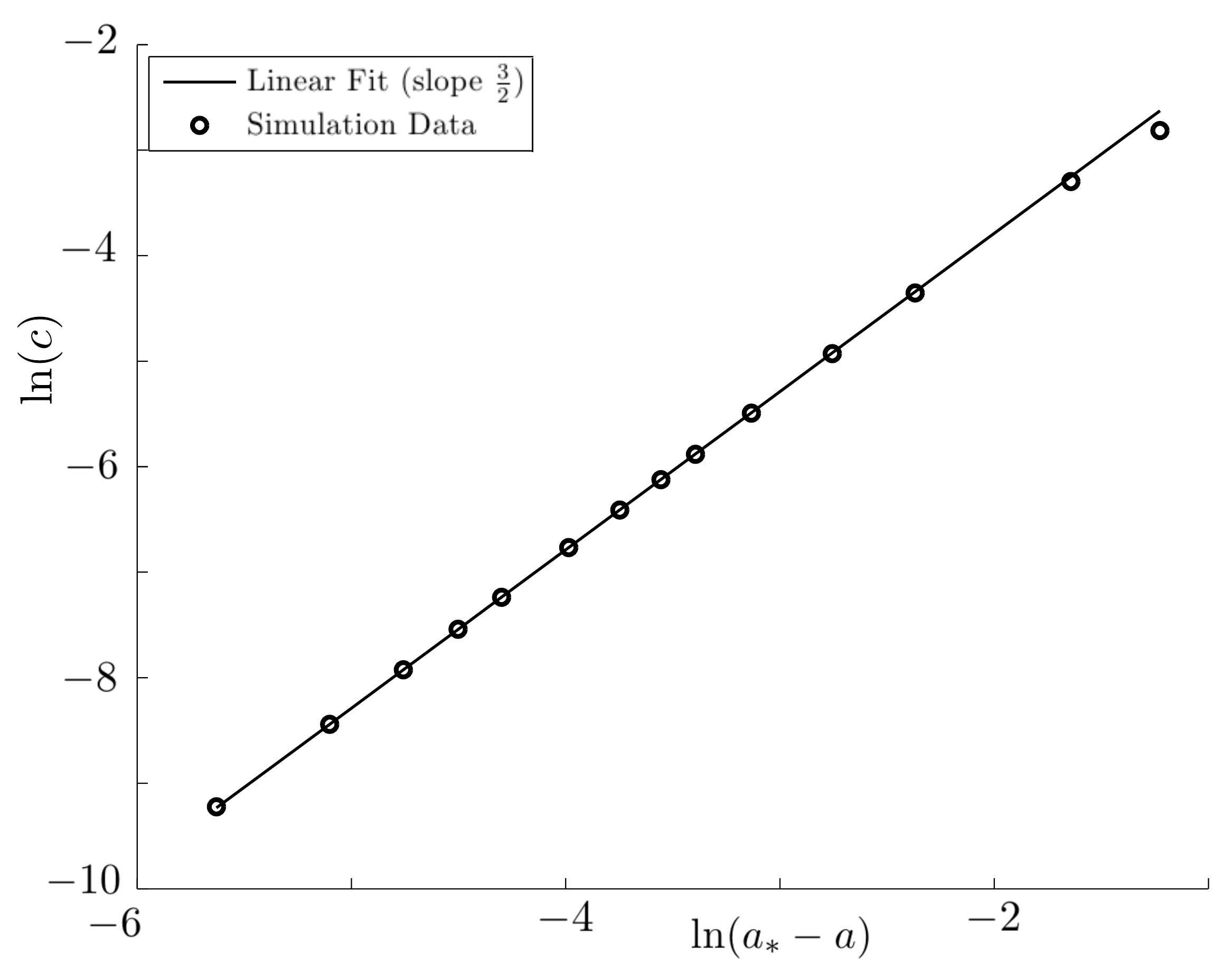}}
%,bb=0 0 398 326
\hspace{.3in}
\subfigure[$\mathcal{K}=\frac{\mathrm{e}^{-x^2}}{\sqrt{\pi}}$, $d=0.1$]{
\label{fig:charKern}
\includegraphics[width=0.45\textwidth]{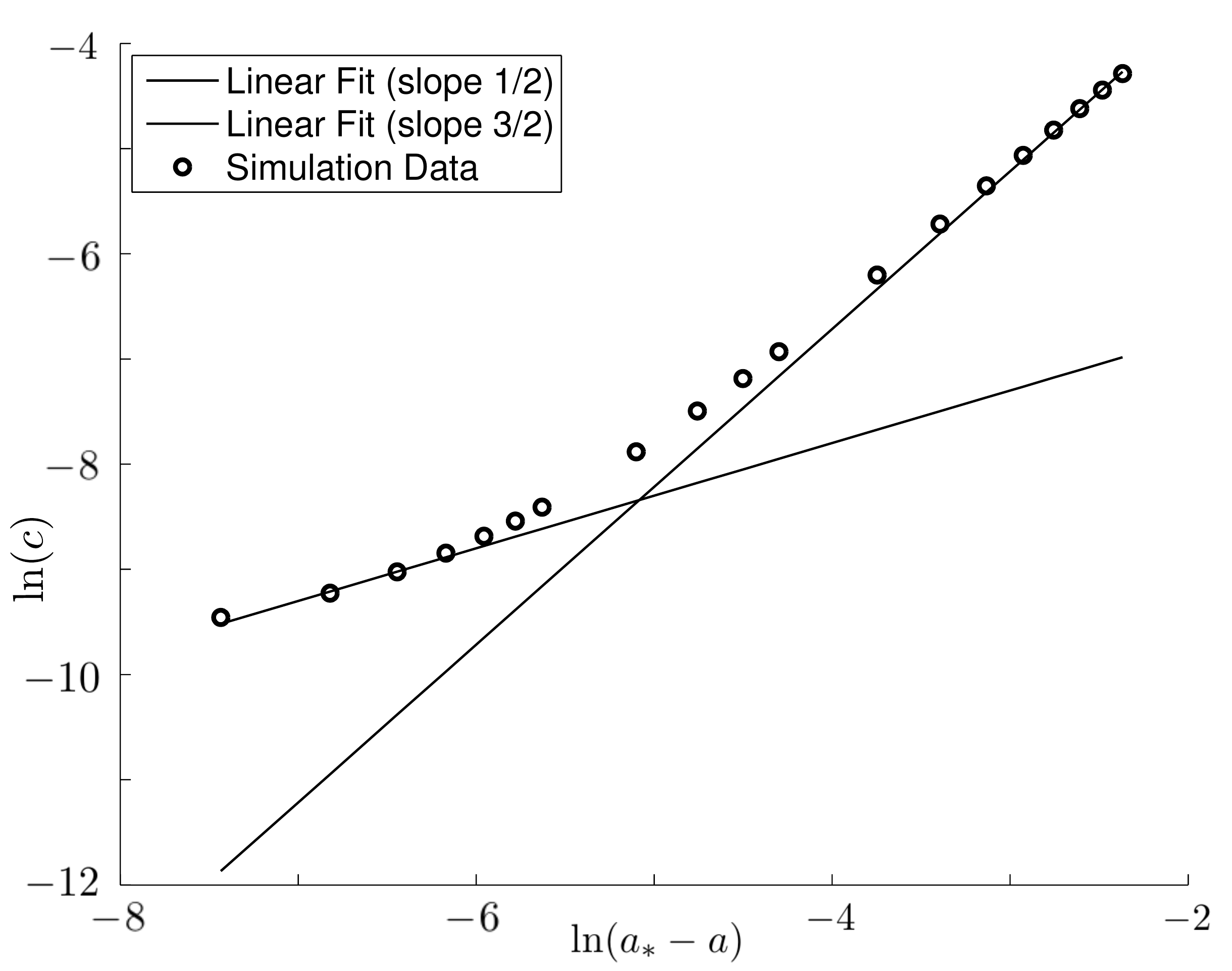}}
%bb=0 0 405 310
\caption{Log-log plots of speeds versus distance to pinning boundary, for various kernels. Plotted are results from \textsc{auto07p} (a) and direct simulations using Fourier transform,  and  best fits for slopes.}
\label{fig:numerics2}
\end{figure}

Figure \ref{fig:numerics2}(a) shows the equivalent results for the first-order, asymmetric kernel. Again, measured scaling law and prefactor are in excellent agreement with the predictions from Theorem \ref{thm:1D}, for a wide range of $\mu$.  
% 
% Throughout, $k_m$ and $k_p$ are measured and predicted scaling factors coefficients in $\eps\sim k(a_\ast-a)^s$, respectively, while $s_m$ and $s_p$ are measured and predicted scaling exponents. Measured values are obtained from linear fits in the log-log plots. 

% To calculate $C_0$ from \eqref{eq:kc} we use an approximate numerical technique; starting from a largely negative point on the cubic $g(u;a,d)$ we integrate backwards in time to find the unique $\gamma_0$ that converges to $g(u;a,d)$ at $\pm\infty$.

\paragraph{Universality: smooth kernels.}

Beyond \textsc{auto07p}, we found that direct simulations produced the most robust results. We therefore explored other kernels (that do not necessarily lead to an ODE formulation) using direct simulations. Specifically, we explored 
\[
\mathcal{K}(x)=\frac{1}{2\pi}\int_\R\frac{1}{1+\ell^2+\ell^4}\rme^{i\ell x}\rmd \ell,\qquad \mathcal{K}(x)=\frac{1}{2}\chi_{[-1,1]}(x),\qquad \mathcal{K}(x)=\frac{\mathrm{e}^{-|x|^2}}{\sqrt{\pi}}.
\]
The results are plotted in Figure \ref{fig:numerics2}(b)--(d). In all cases, we confirmed the power law $3/2$ with reasonable accuracy. The prefactor is not immediately available from heuristics in these cases. The characteristic function gave the poorest match, with a cross-over from a $3/2$- to a $1/2$-scaling law. There are at least two explanations for the somewhat poorer match in this case. First, one expects to see scalings with exponent $1/2$ at small speeds since discretization effects effectively produces a discrete ``lattice ODE'', for which one expects scaling laws of $1/2$ for reasons mentioned in the introduction. On the other hand, one can hope that such discretization effects are smoothed out for smooth kernels, while the evaluation of the characteristic function is most sensitive to grid effects. Measured\footnote{Discretization details: We used $2^{13}$ Fourier modes,  $dt=0.17$, $|x|\leq 30$; best-fit slopes using 4th to 8th smallest data points.}
 and predicted slopes $s_{m/p}$ and $k_{m/p}$ are presented in Table \ref{table:data2}.
%(a) $k_m=0.3892$, $k_p=0.4218$, $s_m=1.492$, $s_p=1.5$ (b) $s_m=1.501$, $s_p=1.5$ (c) $s_m=1.396$, $s_p=1.5$ (d) $s_m=1.498$, $s_p=1.5$\footnote{Discretization details: We used $2^{13}$ Fourier modes,  $dt=0.17$, $|x|\leq 30$; best-fit slopes using 4th to 8th smallest data points.}
\begin{table}[!h] 
\begin{center}
\begin{tabular}{|c||c|c|c|c|}
\hline  & $\mathcal{K}(x)$ & $d$ & $s_m$ & $s_p$   \\\hline
Figure \ref{fig:numerics2}(b) & $(1-\partial_x^2 + \partial_4^4)^{-1}$ & $1/10$  & $1.501$ & $1.5$ \\\hline
Figure \ref{fig:numerics2}(c) & $\chi_{[-1,1]}/2$ & $1/10$ & $1.396$ & $1.5$ \\ \hline
Figure \ref{fig:numerics2}(d) & $\mathrm{e}^{-x^2}/\sqrt{\pi}$ & $1/10$ & $1.498$ & $1.5$\\\hline
\end{tabular}
\end{center}
\caption{Measured and predicted slopes $s_{m/p}$ from Figure \ref{fig:numerics2}(b)-(d).}\label{table:data2} 
\end{table}

\paragraph{Scaling exponent versus kernel smoothness.}

\begin{figure}
\centering
\includegraphics[width=0.5\textwidth]{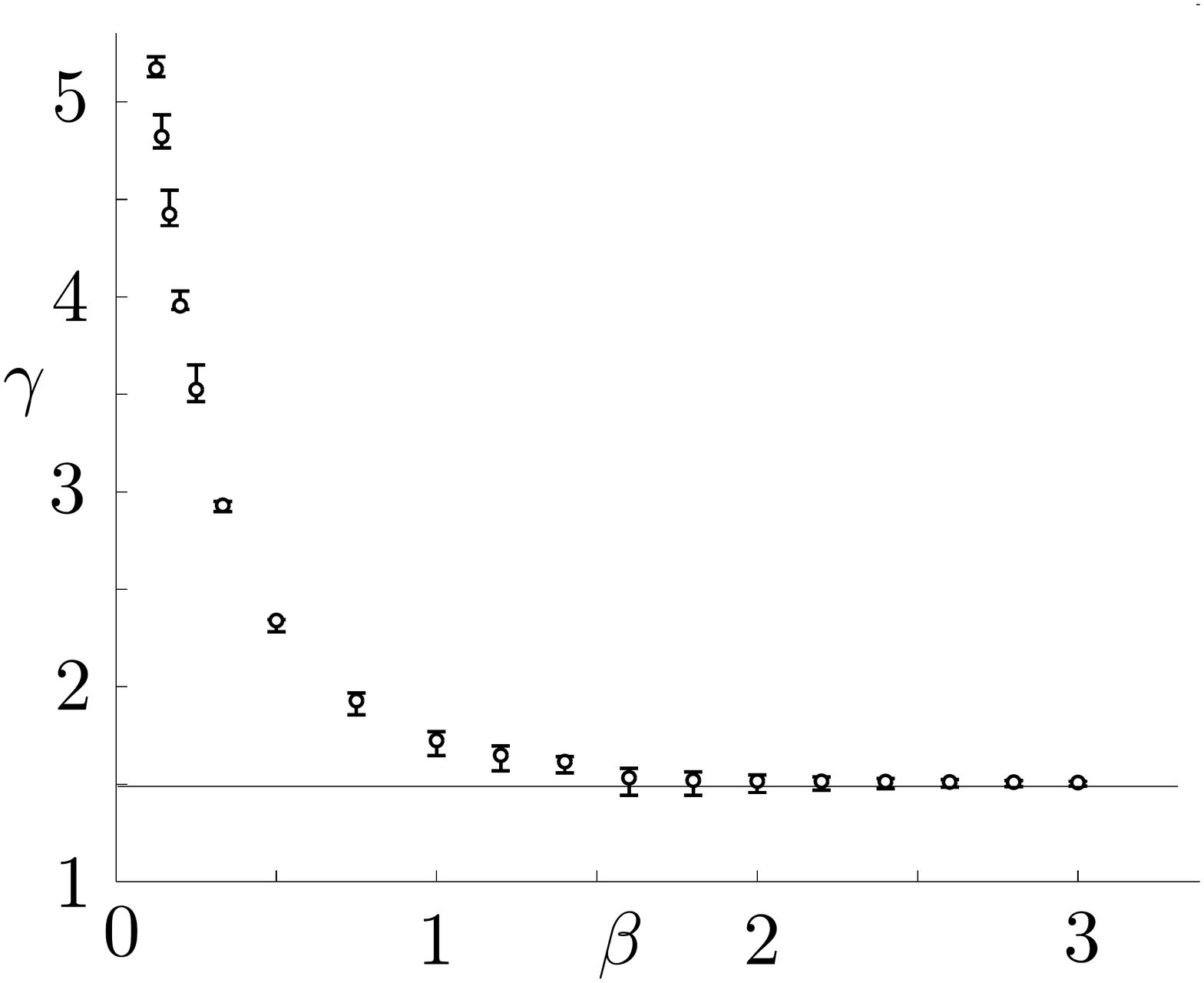}
\caption{Power law exponents $\gamma$ from the relation  $c \sim \mu^\gamma $ for kernels with symbol $\left(1+\ell^2\right)^{-\beta/2}$}
\label{fig:beta}
\end{figure}

Inspired by the prevalence of the scaling exponent $3/2$ in smooth kernels, as opposed to the exponent $1/2$ for discrete, Dirac-delta coupling, we investigated families of kernels with varying degrees of smoothness, 
\[
\widehat{\mathcal{K}}(\ell)=(1+\ell^2)^{-\beta/2}.
\]
Kernels are bounded for $\beta>1$, with a log-singularity at $\beta=1$ and power-law singularities for $\beta<1$. We measured front speeds in direct simulations, using forward Euler in time and Fourier transform in space\footnote{Discretization details: We used $2^{13}$ Fourier modes,  $dt=0.17$, $|x|\leq 30$; best-fit slopes using 4th to 8th smallest data points.}. We found good matches to the speed asymptotics using simple power law relations $c\sim \mu^\gamma$. Figure \ref{fig:beta} summarizes our results. We plotted best fits for scaling exponents as circles and added error bars using maximal and minimal asymptotic slopes. The result is a function $\gamma(\beta)$, which is constant, $\gamma(\beta)\equiv 3/2$ for $\beta\geq 2$. For $\beta<2$, $\gamma$ is monotonically increasing with decreasing $\beta$. 
We also measured scaling exponents for smooth perturbations of these kernels, $(1+\ell^2)^{-\beta/2}+\ell^2\rme^{-\ell^2}$, and $\frac{1}{2}((1+\ell^2)^{-\beta/2}+\rme^{-\ell^2})$. In both cases, the scaling exponent is unchanged by the addition of a smooth part, confirming the observation that scaling laws are primarily determined by the smoothness of the kernel.

\section{Discussion}\label{sec:dis}

We summarize our results and comment on extensions and major open questions. 

\paragraph{Summary of results.} We studied pinning and unpinning in nonlocal problems, comparing pinning regions and unpinning speed asymptotics with discrete and translation invariant problems. For smooth nonlinearities and sufficiently smooth kernels, we found universal scalings for the cuspoidal opening of pinning regions, with width $\sim (d-d*)^2$, where $d_*$ denotes a critical coupling strength; see Theorem  \ref{thm:pinningregion}. Near the pinning boundary, speeds increase with a characteristic power law $\mu^{3/2}$ in most cases. Technically, we showed how geometric desingularization in ODEs can elucidate this scaling for kernels with rational Fourier transform. 

\paragraph{Beyond first and second order kernels.} Our results generalize, at least conceptually, to general rational kernels. In those cases, the fold structure in the fast flow is preserved, while the dynamics on the slow manifold could be more complicated. One still expects similar scalings in the absence of additional degeneracies, such as vanishing of the Melnikov-type derivatives with respect to $a$ in \eqref{eq:k1}, say. It would be very valuable to adapt the geometric desingularization techniques to genuinely nonlinear problems, that is, to problems where the symbol of the kernel is not rational. This would certainly require making the formal asymptotics in \cite{formalbook} rigorous, without using the dynamical systems techniques referenced here \cite{krupa}. Some ideas on how to approach such questions in nonlocal problems can 
be found in \cite{fhnnonlocal}. 

\paragraph{Higher space dimensions.} A straightforward generalization are higher-dimensional systems, with localized kernel $\mathcal{K}(x)$, $x\in\R^n$. Because of translation symmetry, one readily finds one-dimensional profiles $u(e\cdot x - ct)$, for any unit direction $e\in\R^n$. The resulting equation is a one-dimensional non-local problem with ``effective kernel'' $\mathcal{K}_e:=\int_{x\perp e}\mathcal{K}(x)$. For anisotropic kernels, the resulting pinning region may well depend on the direction of propagation, an effect which is well known to lead to faceting of curved interfaces. For discrete, lattice systems, much progress has been achieved more recently by \cite{HP}. It would be interesting to study such problems in the nonlocal context presented here. 

\paragraph{Localization and smoothness.} We have emphasized throughout regularity of the kernel, restricting ourselves to strong, exponential localization. For kernels resulting from $(-\partial_{xx})^{-\beta/2}$, the additional effect of long-range coupling may influence scaling laws. Results in \cite{fraclapl} indicate that speeds would still be finite in this case, but asymptotics are unknown. 

\paragraph{Geometric desingularization.} We relied on the analysis of the simplest singularity of a slow manifold in slow-fast systems to obtain our scaling laws. We also pointed to a slightly more difficult situation, the slow passage through an inflection point, arising in a critical case. It would be very interesting to study the family of problems that arise in these situations more systematically and rigorously using the geometric techniques from \cite{SW}, say. One can, for instance, envision folded nodes and canards in situations where the convolution kernel samples against a nonlinear function of $u$. 

\paragraph{Singular kernels.} Form a theoretical perspective, the most challenging question may be to find and proof power law asymptotics for large classes of singular kernels. Even numerically, such studies are challenging. We are, at this point, not aware of a prediction for the shape of the curve $\gamma(\beta)$ in Figure \ref{fig:beta}. Simple scaling, as used in the passage through the fold or the inflection point, apparently fails to predict correct power laws. Our numerical studies suggest that there may be a fairly simple connection between kernel properties, in particular kernel regularity, and unpinning asymptotics.  

\paragraph{Universality.} A more universal approach to unpinning bifurcations might rely on spectral properties of fronts at the pinning boundary. In fact, we started our paper with the simple picture of a saddle-node, caused by an isolated eigenvalue of the linearization crossing the origin in spatially discrete systems. Embedding those into a continuous system \eqref{e:disnonl} simply increases the multiplicity of this simple eigenvalue, which now is part of the essential spectrum, yet an isolated spectral set. In the case of regularizing convolution operators, one can see, using arguments as in \cite{FS}, that the essential spectrum is given by the numerical range of the derivative of pointwise evaluation operators, $g'(u(x))$. At the pinning boundary, this numerical range touches the origin with a quadratic tangency, for $d<d_*$. This quadratic tangency is reminiscent of quadratic tangencies of spectra near onset of instability in pattern-forming systems, although the spectrum there is parameterized by 
wavenumbers in Fourier space, rather than by the physical location $x$, as in our case here. In this sense unpinning may just be understood as a ``diffusive'', essential instability, depending on interactions with nonlinearity and embedded point spectrum in subtle ways.

\appendix

\section{Passage through an inflection point}

In this appendix, we prove some results stated in Section \ref{subsec:prtip} for the system of differential equations
\begin{subequations}
\begin{align}
\frac{\rmd}{\rmd t}w&=-1,\\
\frac{\rmd}{\rmd t}u&=u^3-w.
\end{align}\label{eq:appendix}
\end{subequations}
This system was obtained after some scalings in the study of the slow passage through an inflection point.
\begin{proposition}
For the system \eqref{eq:appendix}, the following results hold.
\begin{enumerate}
\item There exists a unique trajectory  $\gamma_0$ in the $(u, w)$-plane such that $u(t)^3-w(t)\to 0$ for $t\to \pm \infty$. 
%This trajectory can be written in the form  $\gamma_0=\{(h_0(w),w)\,|\,w \in \R\}$, with $h_0$ monotone, and  $h_0(w) \sim w^\frac{1}{3}$ as $w \rightarrow \pm \infty$. 
\item The Cauchy Principal Value $C_0$ of $u$ exists and
$$C_0 := P.V. \int_\R u(t)\rmd t<0.$$
\end{enumerate}
\end{proposition}

\begin{Proof}

We first start by setting $w_{1,-}:=w^{-\frac{1}{3}}$ and $u_{1,-}:=uw^{-\frac{1}{3}}$ for $w>0$ (respectively $w_{1,+}$ and $u_{1,+}$ for $w<0$) such that system \eqref{eq:appendix} is transformed into two systems
\begin{subequations}
\begin{align}
\frac{\rmd}{\rmd t}w_{1,\pm}&=\frac{1}{3}w_{1,\pm}^4,\\
\frac{\rmd}{\rmd t}u_{1,\pm}&=\frac{u_{1,\pm}^3-1}{w_{1,\pm}^2}+\frac{1}{3}u_{1,\pm}w_{1,\pm}^3.
\end{align}\label{eq:appmod}
\end{subequations}
Now, rescaling time such that $w_{1,\pm}^2 \frac{\rmd}{\rmd t}:=\frac{\rmd}{\rmd s}$, we obtain 
\begin{subequations}
\begin{align}
\frac{\rmd}{\rmd s}w_{1,\pm}&=\frac{1}{3}w_{1,\pm}^6,\\
\frac{\rmd}{\rmd s}u_{1,\pm}&=u_{1,\pm}^3-1+\frac{1}{3}u_{1,\pm}w_{1,\pm}^5.
\end{align}\label{eq:app}
\end{subequations}
Note that each system possesses a unique equilibrium given by $(w_{1,\pm},u_{1,\pm})=(0,1)$, corresponding to the "equilibrium" $(w,u)=(\pm \infty,\pm \infty)$ in system \eqref{eq:appendix}. The linearization at $(0,1)$ for each system is given by the Jacobian matrix 
\begin{equation*}
\mathbf{J}=\left(\begin{matrix} 0 & 0\\ 0 & 3  \end{matrix}\right).
\end{equation*}
As a consequence, for any $k\geq 2$, there exists a $\mathcal{C}^k$ center manifold $\mathcal{M}_{1,\pm}$, given locally, as a graph of form
\begin{equation}
\mathcal{M}_{1,\pm}=\left\{ \left(\Psi_\pm(w_{1,\pm}),w_{1,\pm}\right),~w_{1,\pm} \in \mathcal{V}_\pm \right\},
\label{eq:CenterManifoldApp1}
\end{equation}
where $\mathcal{V}_\pm$ is a neighborhood of the origin in $\mathbb{R}^\mp$, and $\Psi_\pm$ is $\mathcal{C}^k$, with Taylor expansion
\begin{equation}
\Psi_\pm(w_{1,\pm})=1-\frac{1}{9}w_{1,\pm}^5+\mathcal{O}\left(w_{1,\pm}^6\right), \text{ as } w_{1,\pm} \longrightarrow 0.
\label{eq:CenterManifoldApp2}
\end{equation}

We also note that $u_{1,-}\frac{\rmd}{\rmd s}u_{1,-}>0$ provided that $u_{1,-}$ is large enough and $w_{1,-}>0$. As a consequence, $u_{1,-}$ stays bounded as we solve backward in time and using Poincar\'e-Bendixson Theorem we obtain the following result and its corollary.
\begin{lemma}
Fix $0<\delta \ll 1$. Any trajectory of system \eqref{eq:app} with initial condition  $(w_{1,-},u_{1,-})=(\delta,u_{1,-}^0)$, $u_{1,-}^0$ arbitrary, converges backward in time to $(0,1)$. 
\end{lemma}
\begin{corollary}
The asymptotics for $s\longrightarrow -\infty$ are given by the center manifold \eqref{eq:CenterManifoldApp1} and \eqref{eq:CenterManifoldApp2}, up to exponential corrections.
\end{corollary}
Finally, we remark that $u_{1,+} \geq 1$ is locally backward invariant since
\begin{equation*}
\frac{\rmd}{\rmd s}u_{1,+}=\frac{1}{3}w_{1,+}^5 < 0, \text{ at } w_{1,+}<0, u_{1,+}=1.
\end{equation*} 
And for $-S \leq s \leq S$, $S>0$ fixed, solving backward in time, we have that solutions exist ($u_{1,\pm}^2$ decreases when large). 
 
%there exists a unique trajectory that connects $(w_{1,-},u_{1,-})=(0,1)$ at $s=-\infty$ and $(w_{1,+},u_{1,+})=(0,1)$ at $s=+\infty$. 
As a consequence, we have the existence of a unique trajectory $(w,u)$ for system \eqref{eq:appendix} with 
\begin{equation*}
\frac{u^3}{w}\longrightarrow 1 \text{ as } t \longrightarrow \pm \infty,
\end{equation*}
and asymptotics
\begin{equation*}
u(t)=-t^{\frac{1}{3}}\left(1+\frac{1}{9}t^{-\frac{5}{3}}+\mathcal{O}\left(t^{-2}\right) \right),
\end{equation*}
as $t \longrightarrow \pm \infty$. This further ensures that the Cauchy principal value $C_0$ of $u(t)$ exists. One easily checks that $C_0<0$ and this concludes the proof of the proposition.

% As $h_0$ is monotonically increasing, we have that
%\begin{equation*}
%\frac{\rmd}{\rmd w}h_0(w)= \frac{1}{\frac{\rmd}{\rmd z}w}\frac{\rmd}{\rmd z}h_0(w)=-\left( h_0(w)^3-w\right) \geq 0.
%\end{equation*}
%As a consequence, provided that the Cauchy Principal Value of $h_0(w)$ exists, we have
%\begin{equation*}
%P.V. \int_\R h_0(w)\rmd w \leq P.V. \int_\R w^{1/3}\rmd w =0.
%\end{equation*}

\end{Proof}

\bibliography{plain}

\begin{thebibliography}{10}


\bibitem{bates}
P. Bates, P. Fife, X. Ren, and X. Wang,
\newblock {\em Traveling Waves in a Convolution Model for Phase Transitions.} 
\newblock Archive for Rational Mechanics and Analysis, 138, pp. 105--136, 1997.

\bibitem{broer}
H. Broer, T. Kaper and M. Krupa.
\newblock {\em Geometric Desingularization of a Cusp Singularity in Slow-Fast Systems with Applications to Zeeman's Examples}.
\newblock J. Dyn. Diff. Equat., 25, pp. 925--958, 2013.

\bibitem{carpio}
A. Carpio and L. Bonilla.
  \newblock {\em Wave Front Depinning Transition in Discrete One-Dimensional Reaction Diffusion Systems.}
 \newblock Physical Review Letters, 86.26, pp. 6034--6037, 2001. 
 
\bibitem{chen} 
X. Chen.
\newblock {\em Existence, uniqueness, and asymptotic stability of traveling waves in nonlocal evolution equations}.
\newblock  Adv. Differential Equations, 2, pp. 125--160, 1997.

\bibitem{fraclapl}
A. Chmaj.
\newblock {\em Existence of traveling waves in the fractional bistable equation}.
\newblock Archiv der Mathematik, 100,  pp. 473--480, 2013.

\bibitem{FS}
G. Faye and A. Scheel.
\newblock {\em Fredholm properties of nonlocal differential operators via spectral flow.}
\newblock Indiana Univ. Math. J., in press, 2013.  

\bibitem{fhnnonlocal} 
G. Faye and A. Scheel.
\newblock {\em Existence of pulses in excitable media with nonlocal coupling.}
\newblock Submitted, 2013.  

\bibitem{fenichel:79}
N. Fenichel.
\newblock{\em Geometric singular perturbation theory for ordinary differential equations.}
\newblock J. Differential Equations, 31, pp. 53--98,  1979.

\bibitem{gils}
S. van Gils, M. Krupa, and P. Szmolyan.
\newblock{\em Asymptotic Expansions Using Blow-Up.}
\newblock  ZAMP, 56.3, pp. 369--397, 2005.

\bibitem{HP}
A. Hoffman and J. Mallet-Paret.
\newblock
\textit{Universality of crystallographic pinning.}
\newblock J. Dynam. Differential Equations, 22, pp. 79--119, 2010.

\bibitem{krupa}
M. Krupa and P. Szmolyan.
\newblock {\em Geometric Singular Perturbation Theory to Nonhyperbolic Points - Fold and Canard Points in Two Dimensions.}
\newblock SIAM J. Math. Anal., 33, pp. 286--314, 2001.

\bibitem{mp1}
J. Mallet-Paret.
\newblock \textit{The global structure of traveling waves in spatially discrete dynamical systems.}
\newblock J. Dynam. Differential Equations, 11, pp. 4--127, 1999.

\bibitem{mprev} 
J. Mallet-Paret.\newblock
\textit{Traveling waves in spatially discrete dynamical systems of diffusive type.} \newblock Dynamical systems, pp. 231--298, 
Lecture Notes in Math., 1822, Springer, Berlin, 2003. 

\bibitem{formalbook}
E.  Mishchenko, Y. Kolesov, A. Kolesov, and N. Rozov.
\newblock {\em Asymptotic Methods in Singularly Perturbed Systems}.
\newblock Monogr. Contemp. Math., Consultants Bureau, New York, 1994.

\bibitem{schvv} 
A. Scheel and E. van Vleck.
\newblock {\em Lattice differential equations embedded into reaction-diffusion systems}.
\newblock Proc. Royal Soc. Edinburgh A, 139A, pp. 193--207, 2009.

\bibitem{SW}
P. Szmolyan and M. Wechselberger,
\newblock {\em Relaxation Oscillations in $\R^3$}.
\newblock J. Differential Equations, 200, pp. 69--104, 2004. 



\end{thebibliography}

\end{document}